\documentclass[11pt]{article}

\usepackage[T1]{fontenc}
\usepackage{lmodern}
\usepackage{amsmath,amssymb,mathtools}
\usepackage{amsthm}
\usepackage{mathrsfs}
\usepackage[colorlinks]{hyperref}
\usepackage[noabbrev]{cleveref}
\usepackage{graphicx}
\usepackage{float,wrapfig,epsfig}
\usepackage[font=small, labelfont=bf]{caption}
\usepackage[font=small]{subcaption}
\usepackage{cancel}
\usepackage{makeidx}
\usepackage{fancyhdr}
\usepackage{booktabs}
\usepackage{outlines}
\usepackage{cite}
\usepackage{bigdelim}
\usepackage{enumitem}
\usepackage{xcolor}
\usepackage{authblk}

\newtheoremstyle{remboldstyle}
{}{}{}{}{\bfseries}{.}{.5em}{{\thmname{#1 }}{\thmnumber{#2}}{\thmnote{ (#3)}}}
\theoremstyle{remboldstyle}
\newtheorem{rembold}{Remark}[section]

\newtheorem{example}{Example}[section]

\newcommand{\bsym}{\boldsymbol}

\numberwithin{equation}{section}
\hypersetup{urlcolor=blue, colorlinks=true}

\def\keywords{\vspace{.5em}
{\textit{Keywords}:\,\relax%
}}

\graphicspath{{figures/}}

\author[$\dagger$]{Oscar P. Bruno} \author[$\dagger$]{Jagabandhu Paul}
\affil[$\dagger$]{Computing and Mathematical Sciences, Caltech,
  Pasadena, CA 91125} \date{}
\begin{document}
\title{Two-dimensional Fourier Continuation and applications}
\maketitle

\begin{abstract}
  This paper presents a ``two-dimensional Fourier Continuation''
  method (2D-FC) for construction of bi-periodic extensions of smooth
  non-periodic functions defined over general two-dimensional smooth
  domains. The approach can be directly generalized to domains of any
  given dimensionality, and even to non-smooth domains, but such
  generalizations are not considered here. The 2D-FC extensions are
  produced in a two-step procedure. In the first step the {\em
    one-dimensional} Fourier Continuation method is applied along a
  discrete set of outward boundary-normal directions to produce, along
  such directions, continuations that vanish outside a narrow interval
  beyond the boundary.  Thus, the first step of the algorithm produces
  ``blending-to-zero along normals'' for the given function values.
  In the second step, the extended function values are evaluated on an
  underlying Cartesian grid by means of an efficient, high-order
  boundary-normal interpolation scheme. A Fourier Continuation
  expansion of the given function can then be obtained by a direct
  application of the two-dimensional Fast Fourier Transform
  (FFT). Algorithms of arbitrarily high order of accuracy can be
  obtained by this method.  The usefulness and performance of the
  proposed two-dimensional Fourier Continuation method are illustrated
  with applications to the Poisson equation and the time-domain wave
  equation within a bounded domain. As part of these examples the
  novel ``Fourier Forwarding'' solver is introduced which, {\em
    propagating plane waves as they would in free space} and relying
  on certain boundary corrections, can solve the time-domain wave
  equation and other hyperbolic partial differential equations {\em
    within general domains} at computing costs that grow {\em
    sublinearly} with the size of the spatial discretization.
\end{abstract}

\section{Introduction\label{sec:introduction}}
This paper presents a ``two-dimensional Fourier Continuation'' method
(2D-FC) for construction of bi-periodic extensions of smooth
non-periodic functions defined over general two-dimensional smooth
domains. The approach can be directly generalized to domains of any
given dimensionality, and even to non-smooth domains, but such
generalizations are not considered here. The usefulness and
performance of the proposed two-dimensional Fourier Continuation
method are illustrated with applications to the Poisson equation and
the time-domain wave equation within a bounded domain. As part of
these examples the novel ``Fourier Forwarding'' solver is introduced
which, {\em propagating plane waves as they would in free space} and
relying on certain boundary corrections, can solve the time-domain
wave equation and other constant-coefficient hyperbolic partial
differential equations {\em within general domains} at computing costs
that grow {\em sublinearly} with the size of the spatial
discretization.

The periodic-extension problem has actively been considered in the
recent literature, in view, in particular, of its applicability to the
solution of various types of Partial Differential Equations
(PDE)~\cite{Huybrechs_2010,Adcock_2014,ALBIN20116248,amlani-bruno-fc-spectral-2016-307,ASKHAM20171,BRUNO20102009,FRYKLUND201857,LYON20103358,STEIN2016252}. The
contributions~\cite{ALBIN20116248,amlani-bruno-fc-spectral-2016-307,BRUNO20102009,LYON20103358},
in particular, utilize the Fourier Continuation (FC) method in one
dimension in conjunction with dimensional splitting for the treatment
of multidimensional PDE problems.  The dimensional splitting is also
used in~\cite{ELGINDY2019372} to produce Fourier extensions to
rectangular domains in two dimensions, where the Fourier Continuation
is effected by separately applying the one-dimensional FC-Gram
method~\cite{ALBIN20116248,amlani-bruno-fc-spectral-2016-307,BRUNO20102009}
first to the columns and then to the rows of a given data matrix of
function values.  The method does assume that the given smooth
function is known on a rectangular region containing the domain for
which the continuation is sought.

The approach to periodic function extension presented
in~\cite{ASKHAM20171,STEIN2016252} is based on the solution of a
high-order PDE, where the extension shares the values and normal
derivatives along the domain boundary.
Reference~\cite{FRYKLUND201857}, in turn, presents a
function-extension method based on use of radial basis functions. In
that approach, overlapping circular partitions, or patches, are placed
along the physical boundary of the domain, and a local extension is
defined on each patch by means of Radial Basis Functions (RBFs).  A
second layer of patches is placed outside the first, on which the
local values are set to vanish. The zero patches are used in
conjunction with a partition of unity function to blend the local
extensions into a global counterpart. The choice of functions used to
build-up the partition of unity determines the regularity of the
extended function.

The 2D-FC extensions proposed in this paper are produced in a two-step
procedure. In the first step the {\em one-dimensional} Fourier
Continuation method~\cite{amlani-bruno-fc-spectral-2016-307} is
applied along a discrete set of outward boundary-normal directions to
produce, along such directions, continuations that vanish outside a
narrow interval beyond the boundary.  Thus, the first step of the
algorithm produces ``blending-to-zero along normals'' for the given
function values.  In the second step, the extended function values are
evaluated on an underlying Cartesian grid by means of an efficient,
high-order boundary-normal interpolation scheme. A Fourier
Continuation expansion of the given function can then be obtained by a
direct application of the two-dimensional FFT algorithm. Algorithms of
arbitrarily high order of accuracy can be obtained by this method.
Since the continuation-along-normals procedure is a fixed cost
operation,  the cost of the method grows only linearly with the
size of the boundary discretization.

As mentioned above, this paper demonstrates the usefulness of the
proposed general-domain 2D-FC technique via applications to both, the
Poisson problem for the Laplace equation and the time-domain wave
equation.  In the Poisson case the 2D-FC method is utilized to obtain
a \emph{particular solution} for a given right hand side; the boundary
conditions are then made to match the prescribed boundary data by
adding a solution of the Laplace equation which is produced by means
of boundary-integral methods. The Fourier Forwarding approach, in
turn, uses the 2D-FC method to solve the spatio-temporal PDE in the
interior of the domain and it then corrects the solution values near
the boundary by means of a classical time-stepping solver. The overall
procedure, which utilizes large time-steps for the interior solver and
small CFL-constrained time-steps for the near-boundary solver, runs in
computing times that grow {\em sublinearly} with the size of the
spatial discretization mesh.

It is interesting to note that the primary continuation device in the
2D-FC method, namely, continuation along normals to the domain
boundary, is a {\em one-dimensional procedure}. This
one-dimensional continuation procedure can be utilized in a
generalization of the method to $n$-dimensional domains with
$n> 2$. This is in contrast with other extension methods mentioned
above. For example, the RBF-based extension
method~\cite{FRYKLUND201857} requires solution of boundary problems of
increasing dimensionality as the spatial dimension grows, which, given
the method's reliance on dense-matrix linear algebra for the
local-extension process, could have a significant impact on computing
costs. Similar comments apply to PDE-based extension methods such
as~\cite{STEIN2016252}.

The proposed 2D-FC algorithm performs favorably in the context of
existing related approaches. Specific comparisons with results
presented in~\cite{FRYKLUND201857} are provided in
\Cref{sec:numerical-examples-poisson} for a Poisson problem considered
in that reference. The recent contribution~\cite{FENG2020109391}, in
turn, presents an FFT-based high-order solver for the Poisson problem
for {\em rectangular} domains, namely, Cartesian products of
one-dimensional intervals in either two- or three-dimensional
space. The present 2D-FC based Poisson solver achieves, {\em for
  general domains}, a similar performance (similar accuracy and
computing time) to that demonstrated in~\cite[Tables 3 and
4]{FENG2020109391} under the Cartesian-domain assumption.

This paper is organized as follows. After a brief review of the 1D-FC
method presented in \Cref{sec:fc_1d}, the proposed 2D-FC method is
introduced in \Cref{sec:four-cont-meth}. The two main applications
considered, namely, solution of the Poisson and Fourier-Forwarding for
the wave equation, are presented in
\Cref{sec:solut-poiss-equat,sec:fc_ff}. Finally our conclusions are
presented in \Cref{sec:conclusions}.

\begin{figure}[h]
  \centering
  \begin{subfigure}[t]{.8\textwidth}
    \includegraphics[width = 1\linewidth, trim = {5mm 0 0mm 0},
    clip]{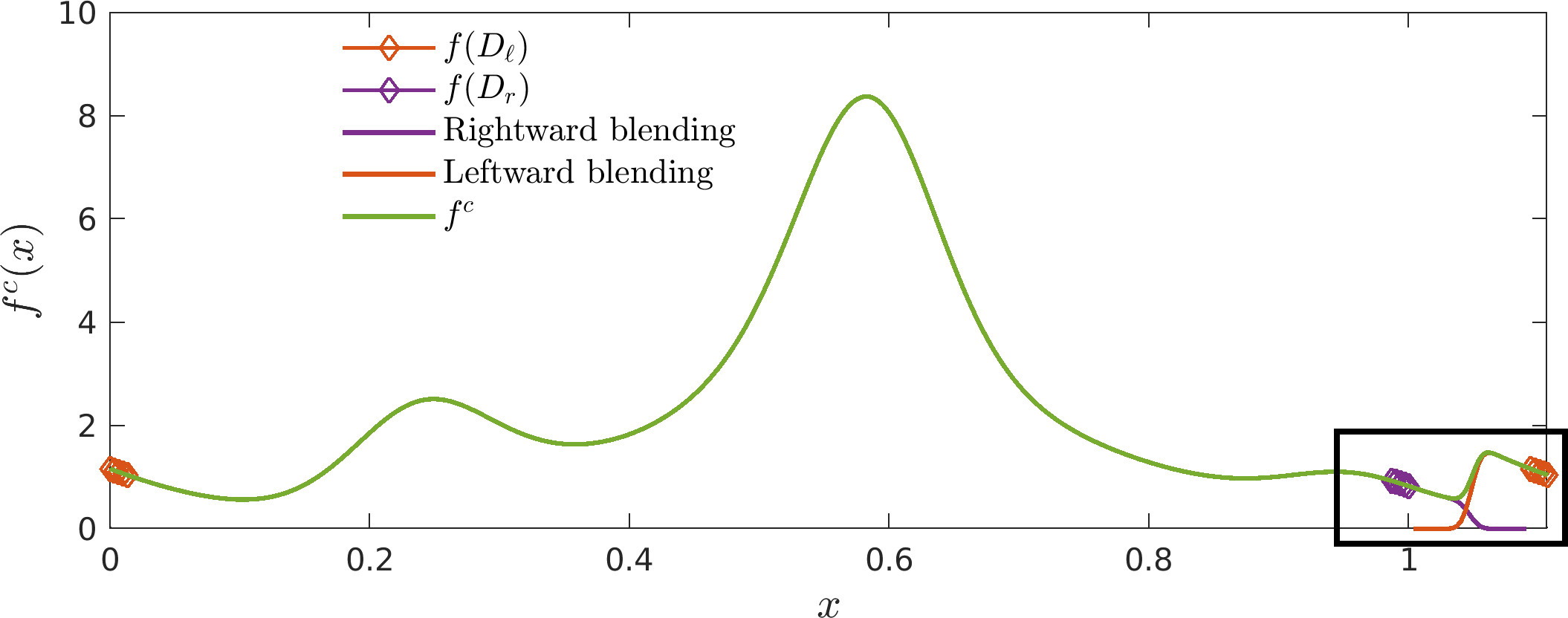}
    \caption{Demonstration of the 1D-FC method. The continuation
      values are computed as the sum of the blended-to-zero rightward and
      leftward extensions.}~\label{fig:fc_1d_values}
  \end{subfigure}
  \begin{subfigure}[t]{0.42\textwidth}
    \includegraphics[width = 1\linewidth, trim = {5mm 0 0mm 0},
    clip]{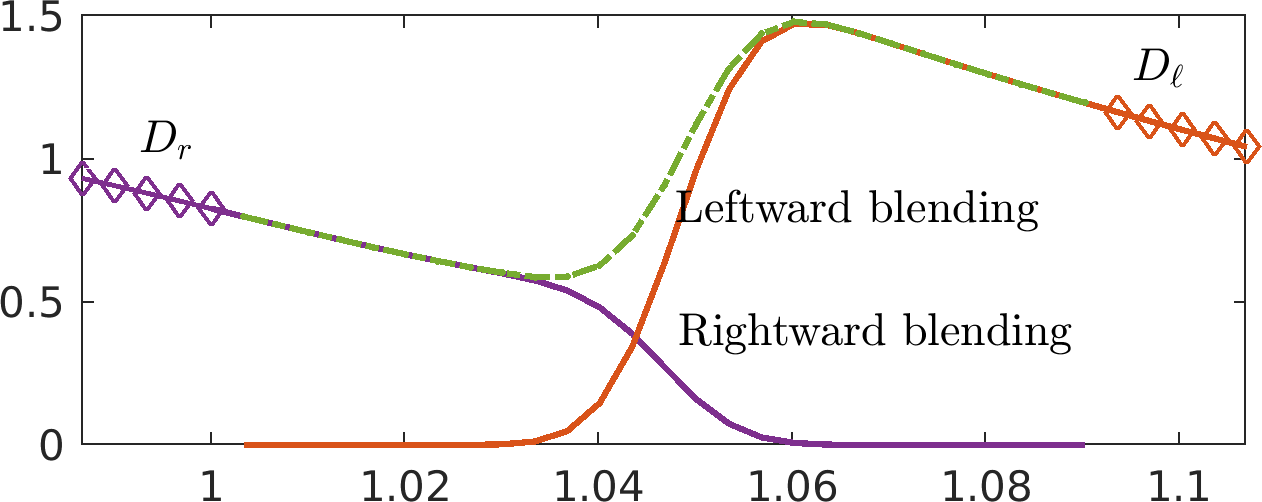}
    \caption{Inset
      of~\Cref{fig:fc_1d_values}. }~\label{fig:fc_1d_blending}
  \end{subfigure}
  \begin{subfigure}[t]{0.4\textwidth}
    \includegraphics[width = 1\linewidth, trim = {5mm 0 0mm 0},
    clip]{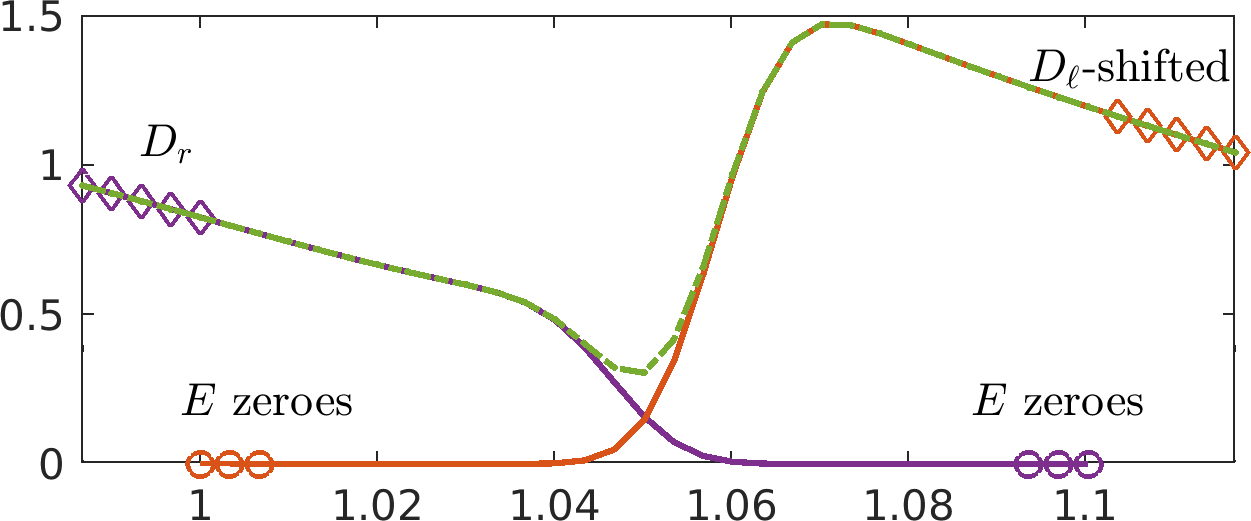}
    \caption{Numbers \(E\) of leftward and rightward extra
      zeroes.}~\label{fig:fc_1d_ezrs}
  \end{subfigure}
  \caption{Illustration of the 1D-FC procedure.
    \Cref{fig:fc_1d_values} depicts the Fourier Continuation of the
    non-periodic function $\phi:[0,1]\to\mathbb{R}$ given by 
    \(\phi(x) = \exp(\sin(5.4\pi x - 2.7\pi) - \cos(2\pi x)) -
    \sin(2.5\pi x) + 1\). \Cref{fig:fc_1d_blending}
    presents a close up of the right continuation region
    \([1-(d-1)k,b]\).  Subsequently \Cref{fig:fc_1d_ezrs} illustrates
    the use of a number \(E\) of extra zeroes in the blending to zero
    process, to yield a continuation mesh containing FFT-friendly
    numbers (products of powers small prime numbers) of point
    values. 
  }~\label{fig:fc_1d}
\end{figure}
\section{Background: 1D ``blending-to-zero'' FC
  algorithm}\label{sec:fc_1d} This section  presents a
brief review of the one-dimensional Fourier Continuation (1D-FC)
method~\cite{amlani-bruno-fc-spectral-2016-307,ALBIN20116248}, with an
emphasis on one of its key components, the \emph{blending-to-zero}
procedure---which is employed in the normal-direction continuation
portion of the proposed two-dimensional Fourier Continuation (2D-FC)
approach presented in~\Cref{sec:four-cont-meth}.

\subsection{1D-FC algorithm: Outline\label{sec:outline}}
Given the vector \(\bsym{\phi}_D = (\phi_0, \dots, \phi_{N - 1})^t\)
of values of a smooth function \(\phi: [0, 1]\rightarrow\mathbb{C}\)
on the equispaced grid \(D = \{x_j = jk: 0\leq j \leq N - 1\}\) of
step-size $k = 1 / (N - 1)$, the 1D-FC
method~\cite{amlani-bruno-fc-spectral-2016-307,ALBIN20116248} of order
\(d\) (with, e.g., $4\leq d\leq 12$) produces, at first, an
(\(N+C\))-dimensional vector $\bsym{\phi}^c$ of discrete continuation
function values (including the \(N\) given function values) over an
extended interval \([0,b], b>1\). To do this, the algorithm utilizes
the \(d\)-dimensional vectors
\(\bsym{\phi}_\ell= (\phi_0, \dots, \phi_{d - 1})^t\) and
\(\bsym{\phi}_r = (\phi_{N-d}, \dots, \phi_{N-1})^t\) of values of the
function $\phi$ on the left and right ``matching-point'' sets
\(D_\ell = \{x_0,\dots,x_{d-1}\}\) and
\(D_r = \{x_{N-d},\dots,x_{N-1}\}\), respectively, each one of which
is contained in a small subinterval of length \((d-1)k\) near the
corresponding endpoint of the containing interval $[0,1]$.  In order
to obtain the \(C\) necessary continuation values, the 1D-FC method
blends $\bsym{\phi}_\ell$ and $\bsym{\phi}_r$ to zero
(see~\Cref{sec:background_fc_1d}), towards the left and right,
respectively, resulting in two zero-blending vectors of length
\(C\). The sum of these two vectors is then utilized as a rightward
discrete continuation to the set
\(D^c=\{x_j = 1+jk: 1\leq j \leq C\}\}\) of points in the interval
\((1,b]\)---as described in~\Cref{sec:fc_1d_algo}. As indicated in
that section, the overall 1D-FC procedure is then completed via an
application of the FFT algorithm to the (\(N+C\))-dimensional vector
$\bsym{\phi}^c$ (cf.~equation~\eqref{eq:fc_1d_cont} below) of
``smoothly-periodic'' discrete continued function values. The
following two subsections describe the blending-to-zero and 1D-FC
approaches, respectively.
\subsection{Blending-to-zero Algorithm}\label{sec:background_fc_1d}
In our description of the order-$d$ blending-to-zero
algorithm~\cite{amlani-bruno-fc-spectral-2016-307} we only present
details for the {\em rightward} blending-to-zero technique, since the
leftward blending-to-zero problem can easily be reduced to the
rightward problem. Thus, given the column vector
\(\bsym{F}_{\mathcal{D}} =(F_0,\dots,F_{d-1})^t \) of values of a
complex-valued function \(F\) on the set
\(\mathcal{D} =\{x_{0},x_1,\dots,x_{d-1}\}\), the rightward
blending-to-zero approach starts by producing a polynomial interpolant
for \(F\) over the interval \([x_{0}, x_{d - 1}]\) relying on the Gram
polynomial basis
\begin{equation}\label{eq:gram_poly}
  G_d = \{g_0(x),g_1(x),\dots,g_{d-1}(x) \}
\end{equation}
for this interval. The functions $g_j(x)$ ($j=0,\dots,d-1$) are the
polynomials with real coefficients that are obtained as the
Gram-Schmidt orthogonalization procedure is applied, in order of
increasing degree, to the polynomials in the set
$\{1,x,x^2,\dots,x^{d-1}\}$, with respect to the discrete
scalar product
\[
  (g, h) = \sum_{j=0}^{d-1}g(x_j)h(x_j).
\]
Discrete values of the
Gram polynomials on the set $\mathcal{D}$ can be computed on the basis of the
\(QR\) factorization~\cite{golub2012matrix}
\begin{equation}\label{eq:pqr_factr}
  \bsym{P} = \bsym{QR}\quad \mbox{of the matrix}\quad
  \bsym{P} = {({x_i} ^ {j - 1})}_{0\leq i, j \leq {d - 1}}. 
\end{equation}
(Note that the \(j ^{\textrm{th}}\) column of \(\bsym{Q}\) contains
the values of the \(j ^ {\textrm{th}}\) Gram polynomial on the set
\(\mathcal{D}\).)  Following~\cite{amlani-bruno-fc-spectral-2016-307} we obtain
the necessary \(QR\) factorization by applying
the stabilized Gram-Schmidt orthogonalization method to the matrix
$\bsym{P}$.

In order to closely approximate each one of the Gram polynomials in
$G_d$, throughout the continuous interval $[0,(d-1)k]$ containing
\(\mathcal{D}\), by corresponding trigonometric polynomials, as described below,
we use a certain ``oversampled matching'' method. According to this
method the polynomials in $G_d$ are oversampled to an equispaced set
of discretization points with stepsize $k/n_{\textrm{os}}$ (containing
\(n_{\textrm{os}}(d-1)+1\) points) where \(n_{\textrm{os}}\) denotes
the oversampling factor, and where the oversampled values are used as
part of a certain Singular Value Decomposition (SVD) matching
procedure described in what follows. Note that the aforementioned
oversampled values on the refined grid
\(\mathcal{D}_{\textrm{os}} \coloneqq \{\tilde{x}_j = j k / n_{os}: 0\leq j
\leq n_{\textrm{os}}(d - 1)\}\)
coincide with the columns of the matrix
\begin{align}\label{eq:q_over}
  \bsym{Q}_{\textrm{os}} &= \bsym{P}_{\textrm{os}} \bsym{R} ^ { - 1}, 
\end{align}
where \(\bsym{P}_{\textrm{os}}\) is the Vandermonde matrix of size
\((n_{\textrm{os}}(d - 1) + 1)\times d\) corresponding to the
oversampled discretization \(\mathcal{D}_{\textrm{os}}\), and where \(\bsym{R}\)
is the upper triangular matrix in equation~\eqref{eq:pqr_factr}.

The aforementioned SVD matching procedure, which is one of the crucial
steps in the FC approach~\cite{amlani-bruno-fc-spectral-2016-307},
produces a band-limited Fourier series of the form
\begin{align}\label{eq:auxi_trig_poly}
  g^c_j(x)
  &= \sum_{m = -M} ^ {M} a_m^j e ^ {\frac{2\pi i m x}{(d + 2C + Z - 1)k}}
\end{align}
for each polynomial \(g_j\in G_d\) ($0\leq j\leq d-1$), where $C$ is
the number of blending-to-zero values to be produced, and $Z$ being the
number of ``zero-matching'' points. The Fourier coefficients are
selected so as to match, in the least square sense, both the
oversampled polynomial values over the interval $[0,(d-1)k]$, and
identically zero values on an equally fine discretization of the
``zero matching'' interval \([(d + C)k, (d + C + Z - 1)k]\) of length
\((Z - 1)k\).
The coefficients in (\ref{eq:auxi_trig_poly}) are taken to equal the
solution $\boldsymbol{a}$ of the minimization problem
\begin{align}\label{eq:min_prob}
  \min_{\bsym{a} = {(a_{ - M}, \dots, a_M)} ^ T}
  \left\|
  \bsym{B} _ {\textrm{os}}\bsym{a} -
  \begin{pmatrix}
    \bsym{q} ^ j_{\textrm{os}}\\
    \boldsymbol{0}
  \end{pmatrix}
  \right\|_{2}, 
\end{align}
where \(\bsym{B}_{\textrm{os}}\) is a matrix whose entries are values
of (\ref{eq:auxi_trig_poly}) at all points in the set
\(\mathcal{D}_{\textrm{os}}\) as well as all the set of
\(k/n_{os}\)-spaced points in the ``zero matching'' interval mentioned
above (which, in particular, contains the endpoints $(d + C)k$ and
$(d + C + Z - 1)k)$. The minimizing Fourier coefficients
\(\boldsymbol{a}\) are then found via an
SVD-based~\cite{golub2012matrix} least-squares approach.\@ Once the
coefficients \(\boldsymbol{a}\) have been obtained, the resulting
Fourier expansions~\eqref{eq:auxi_trig_poly} are used to produce a
certain ``continuation matrix''
\(\bsym{A}\in \mathbb{C}^{C\times d}\), whose columns equal the values
of the expression~(\ref{eq:auxi_trig_poly}) at the $C$ (unrefined)
$k$-discretization points in the interval \([dk, (d+C-1)k]\)
(cf.~\Cref{rmk:Cont_ref} below).  The desired vector $\bsym{F}^r$ of
rightward blending-to-zero function values at the \(C\) continuation
points in the interval \([dk, (d+C-1)k]\) is then given by the
expression
\begin{equation}\label{eq:b_to_z}
  \bsym{F}^r =  \bsym{A}\bsym{Q}^{T} \bsym{F}.
\end{equation}
\subsection{1D-FC Algorithm}\label{sec:fc_1d_algo}
As outlined in \cref{sec:outline}, the 1D-FC algorithm requires use of
a certain rightward (resp.~leftward) blending-to-zero vector
\(\bsym{\phi}_r^r\) (resp.~\(\bsym{\phi}_\ell^\ell\)) for a given
matching-point vector \(\bsym{\phi}_r\)
(resp.~\(\bsym{\phi}_\ell\)). In view of \cref{eq:b_to_z} we define
\(\bsym{\phi}_r^r = \bsym{A}\bsym{Q}^{T} \bsym{\phi}_r\). To obtain
the leftward extension \(\bsym{\phi}_\ell^\ell\), in turn, we first
introduce the ``order reversion'' matrix
${R}^e\in\mathbb{C}^{e\times e}$ ($e\in\mathbb{N}$) by
\[
  {R}^e(g_0,g_1\dots,g_{e-2},g_{e-1})^t =
  (g_{e-1},g_{e-2},\dots,g_1,g_0)^t,
\]
and we then define
\(\bsym{\phi}_\ell^\ell = {R}^C \bsym{A}\bsym{Q}^{T}
{R}^d\bsym{\phi}_\ell\).  A vector \(\bsym\phi^c\) containing both the
$N$ given values in the vector
\(\bsym{\phi}=(\phi_0,\phi_1,\dots,\phi_{N-1})^t\) as well as the
\(C\) ``continuation'' function values is constructed by appending the
sum $\bsym{\phi}_\ell^\ell+\bsym{\phi}_r^r$ at the end of the vector
\(\bsym{\phi}\), so that we obtain
\begin{equation}\label{eq:fc_1d_cont}
  \bsym{\phi}^c_j = 
  \begin{cases}
    \bsym{\phi}_j & \mbox{for } 0\leq j\leq N-1\\
    \left( \bsym{\phi}_\ell^\ell + \bsym{\phi}_r^r\right)_{(j-N)} & \mbox{for
    } N\leq j\leq N+C-1.
  \end{cases}
\end{equation}

Following the various stages of the construction of the vector
\(\bsym\phi^c\) it is easy to check that, up to numerical error, this
vector contains point values of a smoothly periodic function defined
over the interval $[0,b]$. An application of the FFT algorithm to this
vector therefore provides the desired continuation function in the
form of a trigonometric polynomial,
\begin{align}\label{eq:fc_1d_srs}
  \phi^c(x) &= \sum_{\ell=-(N+C)/2}^{(N+C)/2}
              \hat\phi^c_\ell e^{\frac{2\pi i\ell x}{b}},
\end{align}
which closely approximates $\phi$ in the interval $[0,1]$. In fact, as
demonstrated in previous publications
(including~\cite{ALBIN20116248,amlani-bruno-fc-spectral-2016-307}),
for sufficiently smooth functions $\phi$, the corresponding 1D-FC
approximants converge to $\phi$ with order $\mathcal{O}(k^d)$---so
that, as expected, the number $d$ of points used in the
blending-to-zero procedures determines the convergence rate of the
algorithm. (As discussed in these publications, further, in view of
its spectral character the 1D-FC approach enjoys excellent dispersion
characteristics as well.) The 2D-FC algorithm introduced in the
following section also relies on the one-dimensional blending-to-zero
procedure described in \cref{sec:background_fc_1d}, and its
convergence in that case is once again of the order
$\mathcal{O}(k^d)$.

It is important to note that, for a given order $d$, the matrices
\(\bsym{A}\) and \(\bsym{Q}\) can be computed once and permanently
stored on disc for use whenever an application of the blending-to-zero
algorithm is required---as these matrices do not depend on the point
spacing $k$. A graphical demonstration of various elements of the
1D-FC procedure is presented in~\Cref{fig:fc_1d}.

\begin{rembold}[Extra vanishing values]\label{extra_1D}
  The 1D-FC
  implementations~\cite{ALBIN20116248,amlani-bruno-fc-spectral-2016-307}
  allow for an additional number $E\geq 0$ of identically zero
  ``Extra'' function values to be added on a (unrefined)
  $k$-discretization of the interval \([(d+C)k, (d+C+E-1)k]\),
  as illustrated in~\Cref{fig:fc_1d_ezrs}, to obtain a desired overall
  number of discrete function values (including the given function
  values and the continuation values produced) such as e.g., a power
  of two or a product of powers of small prime numbers, for which the
  subsequent application of the fast Fourier transform is particularly
  well suited.  The corresponding use of extra vanishing values for
  the 2D continuation problem is mentioned in Remark~\ref{extra_2D}.
\end{rembold}
\begin{rembold}[Blending to zero on a refined grid]\label{rmk:Cont_ref}
  As indicated above in the present section, the two-dimensional
  Fourier continuation procedure introduced
  in~\Cref{sec:fc-two-dimensions} utilizes the 1D blending-to-zero
  strategy described above in this section to extend a function given
  on a two-dimensional domain $\Omega$ along the normal direction to
  $\Gamma$; the continuation values obtained at all normals
  are then utilized to obtain the continuation function on the
  Cartesian grid by interpolation. As detailed
  in~\Cref{sec:transf-fc-val-2-cartg}, in order to prevent accuracy
  loss in the 2D interpolation step we have found it necessary to use
  1D normal-direction grids finer than the grids inherent in the
  blending-to-zero process itself. To easily provide the necessary
  fine-grid values, a modified fine-grid continuation matrix
  \(\bsym{A}_r\in \mathbb{C}^{C_r\times d}\) is constructed, where
  \(C_r > C\) denotes the number of fine-grid points utilized. The
  modified continuation matrix \(\bsym{A}_r\) can be built on the
  basis of the minimizing coefficients \(\bsym{a}\) in
  (\ref{eq:min_prob}): the corresponding columns of the fine-grid
  continuation matrix \(\bsym{A}_r\) are obtained by
  evaluating~(\ref{eq:auxi_trig_poly}) on the given fine-grid points
  in the interval \(((d-1)k, (d + C - 1)k]\). The necessary
  blending-to-zero function values at the \(C_r\) fine-grid points are
  given by \(\bsym{A}_r\bsym{Q}^T\bsym{\phi}_D\).
\end{rembold}

\section{Two-dimensional Fourier Continuation
  Method}\label{sec:fc-two-dimensions}\label{sec:four-cont-meth} 
This section presents the proposed volumetric Fourier continuation
method on two-dimensional domains \(\Omega\subset \mathbb{R}^2\) with
a smooth boundary \(\Gamma = \overline\Omega\setminus \Omega\), some
elements of which are illustrated in~\Cref{fig:fc_2d}. Let a smooth
function \(f:\overline\Omega\to \mathbb{C}\) be given; we assume that
values of $f$ are known on a certain uniform Cartesian grid within
\(\overline\Omega\) as well as a grid of points on the boundary
\(\Gamma\). The 2D-FC algorithm first produces one-dimensional
``blending-to-zero'' values for the function $f$ {\em along directions
  normal to the boundary} \(\Gamma\), yielding continuation values on
a certain two-dimensional tangential-normal curvilinear grid around
$\Gamma$, as detailed in Sections~\ref{sec:curvilinear}
and~\ref{sec:interp-scheme-bound} and illustrated in~\Cref{fig:fc_2d}.
These continuation values, which are produced on the basis of the
corresponding blending-to-zero procedure presented
in~\cref{sec:background_fc_1d} in the context of the 1D-FC method, are
then interpolated onto a Cartesian grid around the domain boundary, to
produce a two-dimensional blending-to-zero continuation of the
function $f$. The necessary interpolation from the curvilinear grid to
the Cartesian grid is accomplished by first efficiently obtaining the
foot of the normal that passes through a given Cartesian grid point
$\bsym{r} = (x,y)$ exterior to $\Omega$ and near the boundary $\Gamma$
(Section~\ref{sec:Proximity_map}), and then using a local
two-dimensional interpolation procedure to produce the corresponding
continuation value at the point $\bsym{r}$
(Section~\ref{sec:transf-fc-val-2-cartg}). Once the interpolated
values have been obtained throughout the Cartesian mesh around
$\Gamma$, the desired two-dimensional Fourier continuation function
\begin{align}
  \label{eq:fc_2d_srs}
  f^c(x, y)
  &= \sum_{\ell = -N_x / 2 + 1}^{N_x / 2}\sum_{m = -N_y / 2 + 1}^{N_y / 2}
    \hat{f}^c_{\ell, m} e ^ {2\pi i\left(\frac{\ell x}{L_x} + \frac{m y}{L_y}\right)}
\end{align}
(where \(L_x\) and \(L_y\) denotes the period in the \(x\) and \(y\)
directions, respectively) is obtained by means of a 2D FFT.  Following
the algorithmic prescriptions presented in
Sections~\ref{sec:curvilinear}
through~\ref{sec:transf-fc-val-2-cartg}, a summary of the overall
2D-FC approach is presented in Section~\ref{sec:summary-fc-2d}.

\subsection{Two-dimensional tangential-normal curvilinear
  grid\label{sec:curvilinear}} 
\begin{figure}[h]
  \centering
  \includegraphics[width=0.35\linewidth]{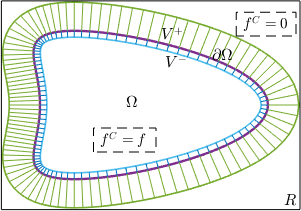}
  \caption{Geometrical constructions underlying the 2D-FC procedure,
    with reference to the various regions defined
    in~\Cref{sec:curvilinear}. }\label{fig:fc_2d}
\end{figure}

The necessary curvilinear grids around
$\Gamma$ can be produced on the basis of a (smooth)
parametrization
\begin{equation}\label{eq:param}
  \bsym{r} = \bsym{q}(\theta)=(x(\theta), y(\theta)), \quad 0\leq \theta\leq 2\pi,
\end{equation}
of the boundary \(\Gamma\). In view of their intended application
(blending to zero along the normal direction in accordance with
\cref{sec:background_fc_1d}), the curvilinear grids are introduced
within interior and exterior strips $V^-$ and $V^+$ (illustrated
in~\Cref{fig:fc_2d}) given by
\begin{equation}\label{eq:int_ext_v}
  \begin{cases}
    V^- = \{\bsym{q}(\theta) - \bsym{n}(\theta) \gamma: 0\leq \theta\leq
    2\pi\mbox{ and }
    0 \leq \gamma \leq (d - 1)k_1  \},\\
    V^+ = \{\bsym{q}(\theta) + \bsym{n}(\theta) \gamma: 0\leq \theta\leq
    2\pi\mbox{ and } 0 \leq \gamma \leq Ck_1\},
  \end{cases}
\end{equation}
where \(\bsym{n}(\theta) = (n_x(\theta), n_y(\theta))\) denotes the unit normal
to the boundary \(\Gamma\) at $\bsym{q}(\theta)$, and where \(d\), $C$
and $k_1 $ denote, respectively, the number of matching points, the
number of blending-to-zero points, and the stepsize used, in the
present application of the 1D blending-to-zero procedure described
in~\Cref{sec:background_fc_1d}.  Using, in addition, the uniform
discretization
\begin{equation}\label{eq:bdry_grid}
  I_{B} = \{\theta_p=p k_2: 0\leq p < B\}, \quad k_2=\frac{2\pi}{B},
\end{equation}
of the interval $[0,2\pi]$, we then construct a curvilinear
two-dimensional discretization
\begin{align}
  \hspace{-4mm} V^-_{B, d} &= \{\bsym{r}_{p, q}: \; \bsym{r}_{p,q}
                             = \bsym{q}(\theta_p) + \bsym{n}(\theta_p)(q - d + 1)k_1;
                             \; 0\leq p < B\mbox{ and }0\leq q\leq d - 1\}\label{eq:vnbj_m} \\
  \hspace{-4mm} V^+_{B, C_r} &= \{\bsym{s}_{p, q}: \; \bsym{s}_{p,q}
                               = \bsym{q}(\theta_p) + \bsym{n}(\theta_p)qk_1 / n_r;
                               \; 0\leq p < B\mbox{ and }0\leq q\leq C_r\},\label{eq:vnbj_p}               
\end{align}
within $V^-$ and $V^+$ respectively, for the given stepsize \(k_1\)
where \(C_r = C n_r\) for certain integer (refinement factor) \(n_r\);
note that the points in \(V^-_{B, d}\) for \(q = d - 1\) and the
points in \(V^+_{B, C_r}\) for \(q = 0\) coincide and that they lie on
\(\Gamma\). Here the constants \(d\), $C$ and
\(n_r\) are independent of $B$.  The continuation function is
constructed so as to vanish at all points
\(\bsym{s}_{p,q}\in V^+_{B, C_r}\) with \(q = C_r\).
Let now \(\mathcal{R} = [a_0, a_1]\times[b_0, b_1]\) denote the smallest closed
rectangle containing \(\Omega\cup V^+\), and consider the equispaced
Cartesian grid of stepsize $h$,
\begin{equation}\label{eq:cart_grid}
  H = \{\bsym{z}_{i, j} = (x_i, y_j): x_i = a_0 + i h;\; y_j = b_0 +
  jh: \; 0\leq i < N_x, 0\leq j < N_y \}
\end{equation}
on \(\mathcal{R}\), where the two-dimensional continuation function
values are to be computed.  We note that the size of the rectangle
\(\mathcal{R}\) along with the strips \(V^-\) and \(V^+\) decrease as
the stepsize $k_2$ is decreased.

\subsection{Computation of FC values on
  \(V^+_{B,C_r}\)}\label{sec:interp-scheme-bound}
\begin{figure}[h]
  \centering \captionsetup[subfigure]{labelformat=empty}
  \begin{subfigure}[t]{0.35\textwidth}
    \includegraphics[width = 1\linewidth, trim = {0 4cm 0 3.2cm},
    clip]{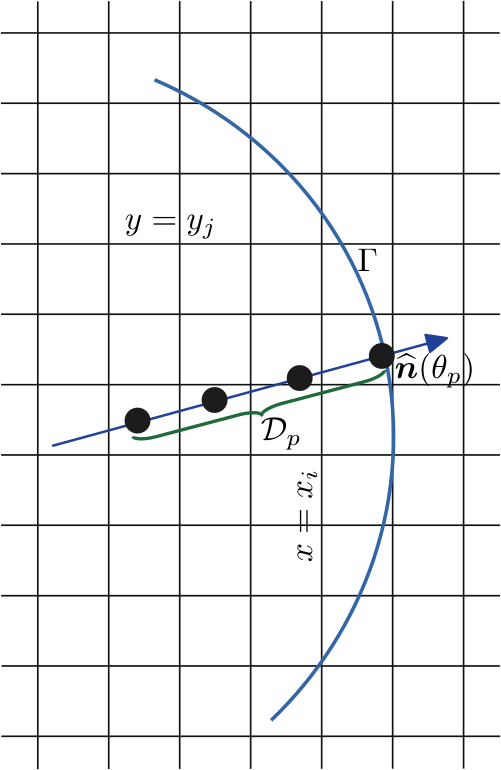} 
  \end{subfigure}
  \hspace{1cm}
  \begin{subfigure}[t]{0.35\textwidth}
    \includegraphics[width = 1\linewidth, trim = {0 4cm 0 3.2cm},
    clip]{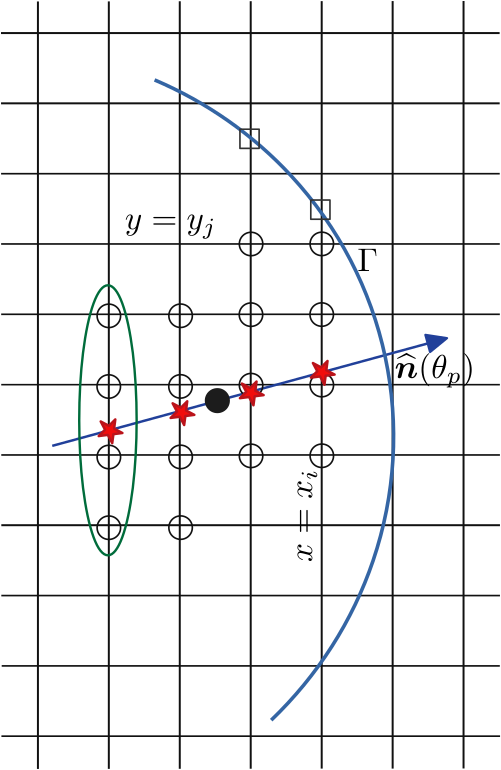}   
  \end{subfigure}
  \caption{Interpolation scheme for evaluation of
    \(\bsym{f}_{\mathcal{D}_p}\) in the case
    \(|n_x(\theta_p)| \geq |n_y(\theta_p)|\).  Left: Black solid
    circles indicate the matching points that define the set
    \(\mathcal{D}_p\). Right: Known function values at the Cartesian
    points (denoted by the empty circles) and, in some cases, at the
    point of intersection (represented by an empty square) of the
    vertical grid lines with \(\Gamma\), are used to interpolate the
    matching function values at the red-star intersection points of
    the normal with the vertical grid lines. The function values at
    the red-points are then used to obtain, by interpolation, the
    function values at the matching
    points. }~\label{fig:bdry-sec-calculation}
\end{figure}

A continuation of the function \(f\) to the exterior of $\Omega$ is
obtained via application of the blending-to-zero procedure presented
in~\Cref{sec:background_fc_1d} (cf.~\Cref{rmk:Cont_ref}) along each
one of the normal directions inherent in the definition of the set
$V^+_{B, C_r}$.  For given \(p\), the \(d\) equidistant points
\(\bsym{s}_{p,q}\in V^-_{B, d}\) \((0\leq q \leq d - 1)\), which are
indicated by the solid circles in \Cref{fig:bdry-sec-calculation},
constitute a set \(\mathcal{D}_p\) of matching points that are used to
effect the blending-to-zero procedure per the prescriptions presented
in~\Cref{sec:background_fc_1d}. To obtain the desired continuation
function values it is necessary to first obtain the vector
\(\bsym{f}_{\mathcal{D}_p}\) of the values of the function $f$ (or
suitable approximations thereof) on the set \(\mathcal{D}_p\). In the
proposed method, the needed function values
\(\bsym{f}_{\mathcal{D}_p}\) are computed on the basis of a two-step
polynomial interpolation scheme, using polynomials of a certain degree
(\(M-1\)), as briefly described in what follows.

With reference to the right image of
\Cref{fig:bdry-sec-calculation}, and considering first the case
\(|n_x(\theta_p)| \geq |n_y(\theta_p)|\), the algorithm initially
interpolates vertically the function values at \(M\) open-circle
Cartesian points selected as indicated in the figure, onto the points
of intersection, shown as red-stars, of the normal line and the
vertical Cartesian grid lines. For intersection (red-star) points
close enough to the boundary, boundary function values at boundary
points shown as squares in the figure, are utilized in the
one-dimensional interpolation process as well.  Once the red-star
function values are acquired, the function value at the matching
solid-black point is effected by interpolation from the \(M\) red-star
point values previously obtained, on the basis of a polynomial of
degree (\(M-1\)).

The case \(|n_x(\theta_p)| < |n_y(\theta_p)|\) is treated similarly,
substituting the initial interpolation along vertical  Cartesian
lines, by interpolation along horizontal Cartesian lines; the
algorithm then proceeds in an entirely analogous fashion. 
\subsection{Proximity map}\label{sec:Proximity_map} 
As described below in Section~\ref{sec:transf-fc-val-2-cartg}, the 2D
FC algorithm interpolates the Fourier continuation values on
\(V^+_{B, C_r}\) onto the Cartesian mesh $H$. The interpolation
algorithm used in that section relies on a certain ``proximity map''
$\mathcal{P}:H\cap V^+\to V^+_{B,C_r}$ which associates a
curvilinear grid point
\(\bsym{s}_{p,q} = \mathcal{P}(\bsym{z}_{i,j}) \in V^+_{B,C_r}\) in the
``proximity'' of each given Cartesian grid point $\bsym{z}_{i,j}$. The
proximity function we use is obtained by first associating to each
curvilinear discretization point $\bsym{s}_{p,q}$ the nearest
Cartesian point, a procedure that results in a set
$\mathcal{P}_0 \subseteq (H\cap V^+)\times V^+_{B,C_r}$ of pairs of
points, one in the Cartesian grid and the other in the curvilinear
grid. (The initial set $\mathcal{P}_0$ can easily be obtained by using
the ``integer part'' \emph{floor} \((\lfloor.\rfloor)\) and the
\emph{ceil} \((\lceil.\rceil)\) operators.) The set $\mathcal{P}_0$ is
then modified by removing multiple associations for a given Cartesian
point, and, if necessary, by adding a ``next-nearest'' curvilinear
neighbor to Cartesian points that previously remained
un-associated. The resulting set $\mathcal{P}$ defines the desired
function. 
\subsection{FC Values on the Cartesian Grid}\label{sec:transf-fc-val-2-cartg} 
Once Fourier continuation values on \(V^+_{B,C_r}\) have been
obtained, per the procedure presented in
\Cref{sec:interp-scheme-bound}, the two-dimensional FC scheme can be
completed by (a)~Interpolation onto the set \(H\cap V^+\) of outer
Cartesian grid points; and, (b)~Subsequent evaluation of the
corresponding Fourier coefficients in equation~\Cref{eq:fc_2d_srs} by
means of an FFT. (Note that since the continuation function $f^c$ is a
smooth function which vanishes outside \(\overline\Omega\cup V^+\),
this function can be viewed as the restriction of a smooth and
bi-periodic function with periodicity rectangle $\mathcal{R}$---whose
Fourier series approximates $f^c$, and therefore $f$, accurately.)
The efficiency of the interpolation scheme is of the utmost importance
in this context---since interpolation to a relatively large set
\(H\cap V^+\) of Cartesian points is necessary. An accurate and
efficient interpolation strategy is obtained by combining two
one-dimensional (local) interpolation procedures based on nearby
normal directions. The first interpolation procedure produces the
parameter value \(\theta\) of the foot of the normal to \(\Gamma\)
passing through a given Cartesian point; the second procedure then
approximates the continuation function value utilizing the parameter
value \(\theta\) just mentioned and the continuation function values
at the points in \(V^+_{B, C_r}\) around the given Cartesian point. A
detailed description of the combined interpolation methodology is
presented in what follows.  Specifically, we describe the strategy we
use to interpolate the continuation function onto each point
\(Q = \bsym{z}_{i, j}\in H\cap V^+\) and, to do this, we first obtain
the foot $\mathcal{F}(Q)$ of this point. Using the proximity map
\(\mathcal{P}\) described in~\Cref{sec:Proximity_map}, the algorithm
utilizes the curvilinear discretization point
$\bsym{s}_{p, q}= \mathcal{P}(\bsym{z}_{i, j})\in V^+_{B,C_r}$ as well
as the corresponding boundary discretization parameter value
\(\theta_p\); according to \cref{eq:vnbj_p}, the point
$\bsym{q}(\theta_p)$ equals the foot of the normal passing through
$\bsym{s}_{p, q}$:
$\bsym{q}(\theta_p) = \mathcal{F}(\bsym{s}_{p, q})$. The algorithm
then seeks approximation of the foot $\mathcal{F}(Q)$ and the
corresponding parameter value \(\theta = \theta^Q\) via a preliminary
interpolation step, as indicated in what follows. The foot
$\mathcal{F}(Q)$ and the corresponding parameter value
\(\theta = \theta^Q\) are then used by the algorithm to produce the
desired interpolated continuation value at \(Q\).

In order to obtain $\mathcal{F}(Q)$ the algorithm uses the \(M\)
boundary parameter values in the set 
\begin{equation}\label{surround}
  S_{\theta_p} =
  \{\theta_{p-K_\ell},\theta_{p-K_\ell+1},\dots,\theta_p,\theta_{p+1},
  \dots,\theta_{p+K_r}\}\subset I_{B}
\end{equation}
around \(\theta_p\) (where $K_r+K_\ell +1 = M$, and where $K_r=K_\ell$
if $M$ is odd and $K_r=K_\ell+1$ if $M$ is even. Parameter
values $\theta_k$ with negative values of $k$, which may arise in
\Cref{surround}, are interpreted by periodicity: $\theta_k =
\theta_{B+k}$).

\begin{figure}[h]
  \centering
  \captionsetup[subfigure]{labelformat=empty}
  \begin{subfigure}[t]{0.4\textwidth} 
    \includegraphics[width = 1\textwidth,trim = {2.6cm 4.5cm 2cm 1.5cm},
     clip]{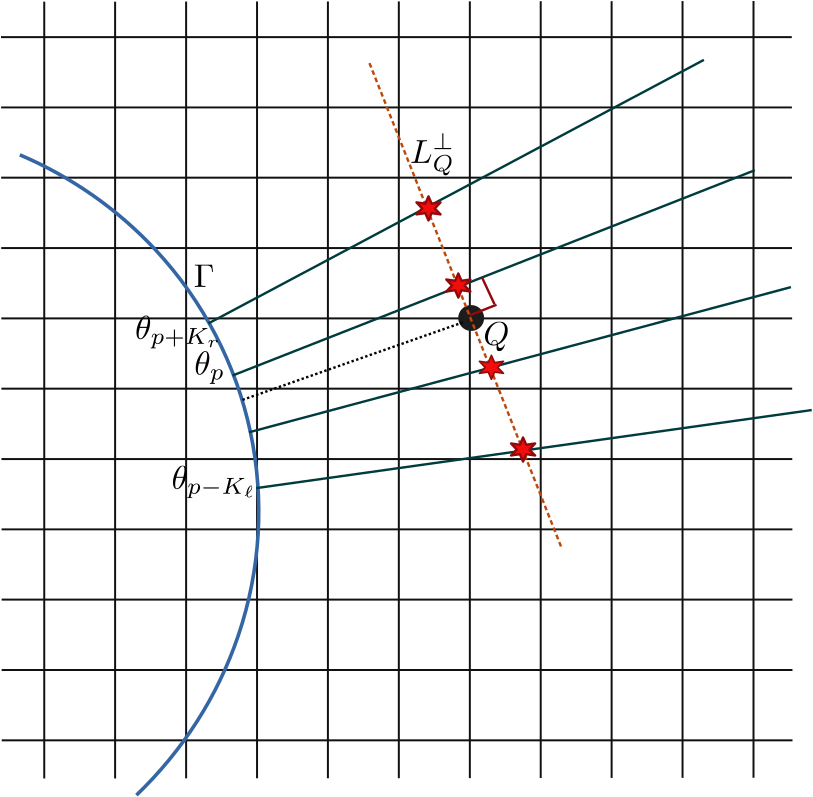}  
  \end{subfigure}
  \hspace{4mm}
  \begin{subfigure}[t]{0.4\textwidth}
    \includegraphics[width = 1\linewidth,trim = {2.6cm 4.5cm 2cm 1.5cm},
     clip]{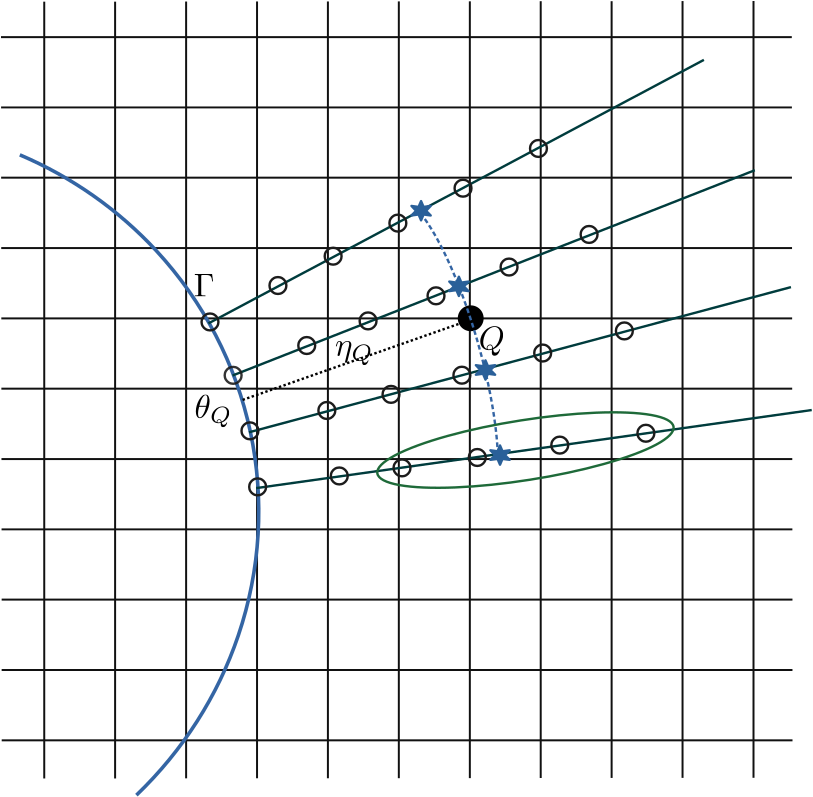} 
  \end{subfigure}
  \caption{Interpolation schemes utilized to obtain the continuation
    function values on \(H\cap V^+\) on the basis of the continuation
    values on \(V^+_{B,C_r}\). Left: Evaluation of the boundary
    parameter value for the foot of the normal line passing through
    \(Q\), depicted as a finely dotted line passing through that
    point. The left image also displays the set of red-star
    interpolation points along the dashed-orange line
    \(L^\perp_Q\). Right: Interpolation of continuation values from
    the curvilinear mesh to a point \(Q\) on the Cartesian
    mesh.}~\label{fig:interpolation_to_cartesian_grid}
\end{figure}

The algorithm then utilizes the line $L^{\perp}_Q$ passing through $Q$
that is orthogonal to the normal vector $\bsym{n}(\theta_p)$ (see left
image in \Cref{fig:interpolation_to_cartesian_grid}), together with
the parametrization \(\bsym{\ell}^{\,\textrm{tr}}_Q(\tau)\) of
$L^\perp_Q$, where the parameter \(\tau\) represents the signed
distance of the points on $L^\perp_Q$ from the point \(Q\). Clearly,
then, \(\bsym{\ell}^{\textrm{tr}}_Q(0) = Q\).  Each point of
intersection of $L^\perp_Q$ with the normals $\bsym{n}(\theta_j)$
($\theta_j\in S_{\theta_p}$), on the other hand, equals
\(\bsym{\ell}^{\textrm{tr}}_{Q}(\tau_j)\) where \(\tau_j\) denotes the
signed distance between \(Q\) and the corresponding intersection
point.  Thus, defining the function \(\theta = \mathcal{T}(\tau)\),
where $\mathcal{T}(\tau)$ gives the parameter value of the foot of the
normal through the point \(\bsym{\ell}^{\textrm{tr}}_{Q}(\tau)\), we
clearly have
\begin{align}\label{eq:param_tau}
  \theta_j = \mathcal{T}(\tau_j); \quad p-K_\ell\leq j\leq p+K_r.
\end{align}
It follows that a 1D interpolation procedure on the function
$\mathcal{T}(\tau)$ can be used to obtain the desired approximation of
the value \(\theta^Q = \mathcal{T}(0)\) of the parameter corresponding to
the foot of the point $Q=\bsym{z}_{i,j}$:
\(\mathcal{F}(Q)  = \bsym{q}(\theta^Q) = \bsym{q}(\mathcal{T}(0))\).

Once we have the corresponding foot parameter value \(\theta^Q\) for
the given Cartesian point \(Q\), the distance \(\eta_Q\) of the point
\(Q\) to the boundary \(\Gamma\) is easily computed. Let
\(S_{\theta^Q}\subset I_{B}\) be a set of \(M\) boundary parameter
values, similar to the set \(S_{\theta_p}\) defined in
(~\ref{surround}), but with \(S_{\theta^Q}\) ``re-centered'' around
\(\theta^Q\). In order to obtain the continuation function value at
the point \(Q\) using the continuation function values on
\(V^+_{B,C_r}\), we employ a local two-dimensional polynomial
interpolation scheme based on the set \(S_{\theta^Q}\) and the
distance \(\eta_Q\), as indicated in what follows and illustrated in
the right image of \Cref{fig:interpolation_to_cartesian_grid}. First,
the continuation values are obtained at each point (marked by blue
asterisks in the figure) at a distance \(\eta_Q\) from the boundary
\(\Gamma\) along the normal grid lines in \(V^+_{B,C_r}\) that
correspond to boundary parameter values in the set
\(S_{\theta^Q}\). Each one of these values is obtained via
one-dimensional interpolation of the continuation function values on
\(V^+_{B,C_r}\) along the corresponding normal grid line in
\(V^+_{B,C_r}\).  The desired continuation value at the point \(Q\),
then, is obtained via a final one-dimensional degree-(\(M-1\))
interpolation step based on the parameter set \(S_{\theta^Q}\) and
values at the ``blue-asterisk'' points just obtained.
Finally, by applying the two-dimensional FFT to the continuation
function values computed above we then obtain the desired Fourier
series expression in (\ref{eq:fc_2d_srs}) for the continuation
function.
\begin{rembold}[Function values on \(\Gamma\)]
  I.~For definiteness, in this paper we have assumed that the boundary
  data is provided in the form of values of the given function---which
  corresponds to the the Dirichlet boundary data in the PDE
  context. But the approach is also applicable in cases for which the
  boundary data is given as the normal derivative of the function,
  (Neumann boundary data), or even a combination of function and
  normal derivative values (Robin data) by relying on a slightly
  modified blending-to-zero procedure of the type presented
  in~\cite{amlani-bruno-fc-spectral-2016-307}.  II.~If no boundary
  data available, the two-dimensional Fourier continuation method can
  still be utilized on the basis of interior data only, albeit with a
  certain reduction in accuracy near the boundary.\label{rmk:fc_2d_bc}
\end{rembold}
\begin{rembold}[Extra vanishing values in 2D]\label{extra_2D}
  As in the 1D case, prior to the FFT procedure the grid $H$ can be
  enlarged, with vanishing function values assigned to the added
  discretization points to obtain a discretization containing a number
  of discretization points equal to a power of two (or a product of
  powers of small prime numbers) along each Cartesian direction, which
  leads to specially fast evaluations by means of the fast Fourier
  transform.
\end{rembold}
\begin{figure}[h]
  \centering
  \begin{subfigure}[b]{0.49\textwidth}
    \includegraphics[width=1\linewidth]{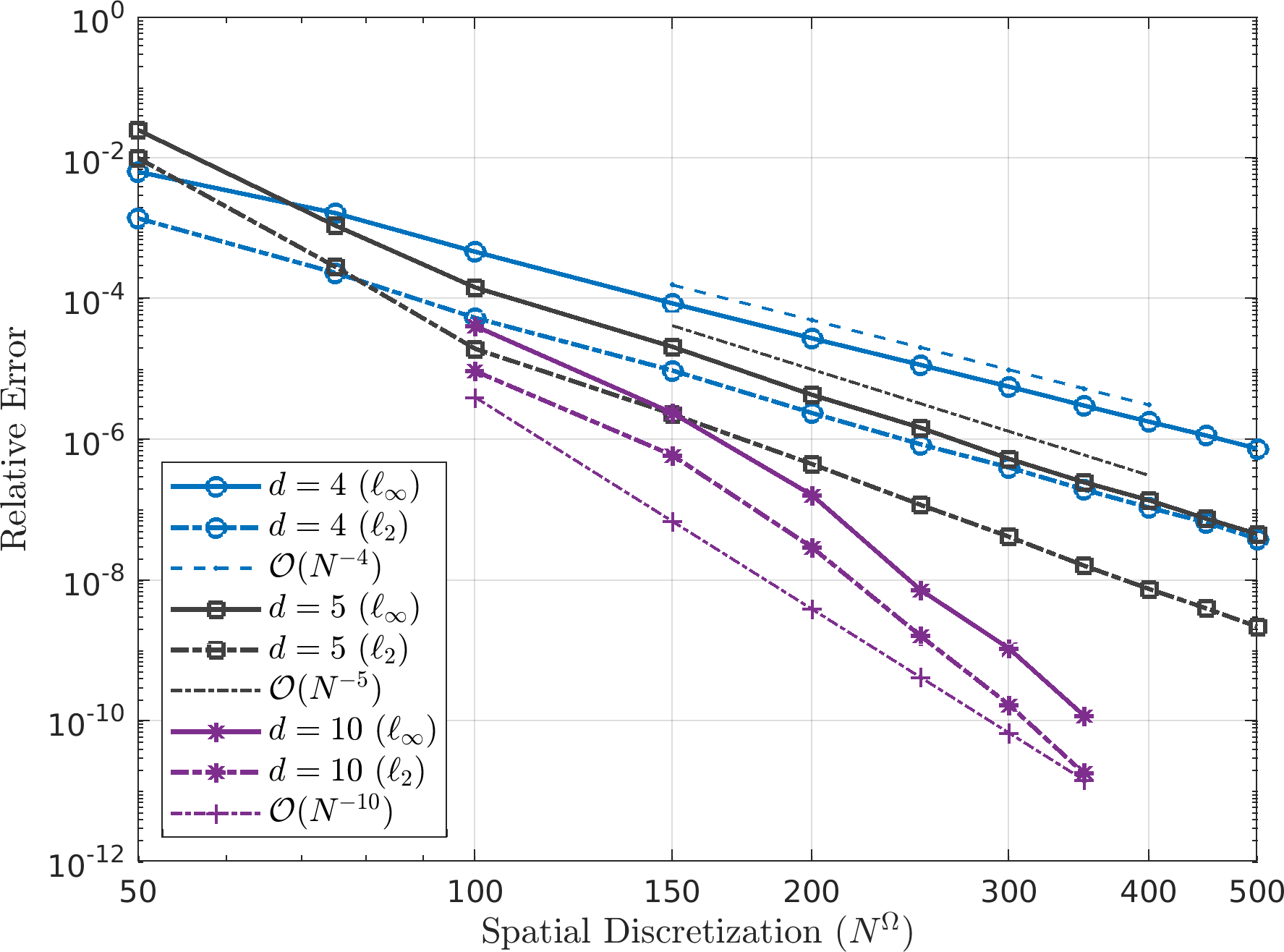} 
  \end{subfigure}
  \begin{subfigure}[b]{0.49\textwidth} 
    \includegraphics[width = 1\linewidth]{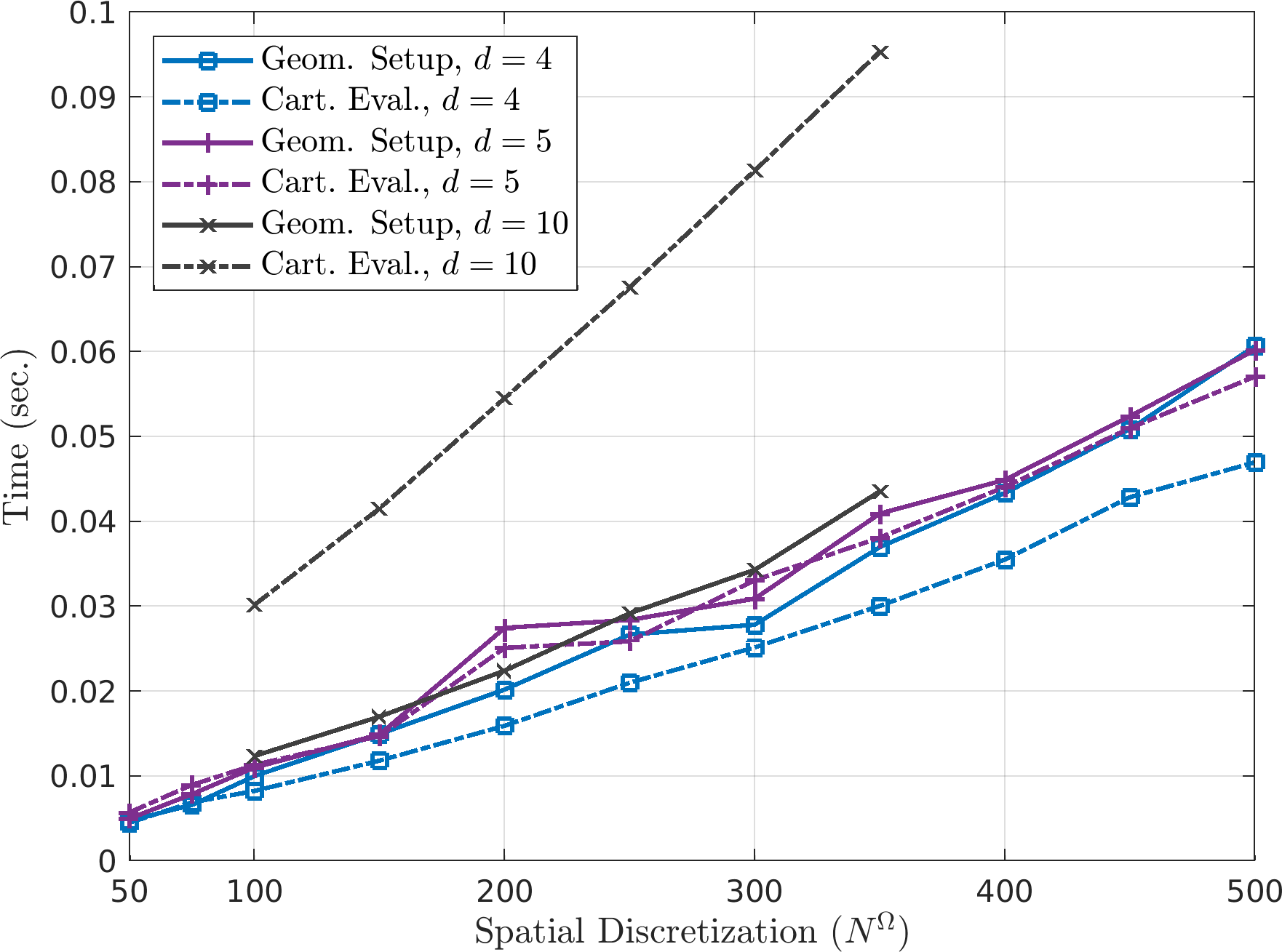}   
  \end{subfigure}
  \caption{Numerical errors in \(\log\)-\(\log\) scale (left graph)
    and computing times required (right graph) in the 2D-FC
    approximation of the function \(f\) in~\eqref{eq:fc_example_1} in
    the setting of \Cref{exm:fc_2d_disc_1}. The interpolating
    polynomial degree \(M=d+3\) was used in all cases. The integer
    \(N^\Omega\) equals the number of spatial grid points used over
    the diameter of the disc. Times reported correspond to averages
    over \(10\) runs.  }~\label{fig:fc_2d_err_time_4_5_10}
\end{figure}
\subsection{Summary of the 2D-FC procedure}\label{sec:summary-fc-2d}
This section presents a summary of the 2D-FC procedure described in
the Sections~\ref{sec:curvilinear}
through~\ref{sec:transf-fc-val-2-cartg} for a function $f$ given on a
uniform Cartesian grid \(H\cap \Omega\) within the domain of
definition $\Omega$, where $H$ is a Cartesian mesh over the rectangle
\(\mathcal{R}\) containing both \(\Omega\) and the near-boundary outer
region \(V^+\); see~\Cref{sec:curvilinear} and, in particular,
\Cref{fig:fc_2d}.  The construction of the continuation function
\(f^c\) for the given function \(f\) relies on use of three main
parameters associated with the 1D blending-to-zero approach presented
in \Cref{sec:fc_1d}, namely \(d\) (number of points in the set
\(\mathcal{D}\)), \(C\) (number of unrefined discrete continuation
points), and \(n_{os}\) (oversampling factor for the 1D
blending-to-zero FC procedure), together with the parameters \(n_r\)
(refinement factor for the discrete continuation
points,~\Cref{rmk:Cont_ref} and~\Cref{sec:curvilinear}), and, \(M-1\)
(degree of the interpolating polynomials,
\Cref{sec:interp-scheme-bound,sec:transf-fc-val-2-cartg}).
Additionally, the 2D-FC procedure utilizes the precomputed matrices
\(\bsym{A}_r\) and \(\bsym{Q}\), which, with reference to
\Cref{sec:background_fc_1d}, are obtained as per the description
provided in \Cref{rmk:Cont_ref}.

Using the aforementioned parameters and matrices, the algorithm
proceeds in two main steps, namely step (a) A {\em ``Geometrical
  Setup''} precomputation procedure (comprising points 1.~through
4.~below); and step (b) A {\em ``Cartesian Evaluation''} procedure
(comprising points 5.~through 7.~below). Part~(a) only concerns
geometry, and, for a given domain and configuration, could be produced
once and used for evaluation of Fourier continuations for many
different functions, as is often necessary in practice. The full 2D-FC
algorithm thus proceeds according to the following seven-steps:
\begin{enumerate}[noitemsep]
\item Discretize the boundary \(\Gamma\) using a smooth
  parametrization
  \(\{\bsym{q}(\theta)=(x(\theta), y(\theta)): 0\leq \theta\leq 2\pi\}\) of
  \(\Gamma\) and the uniform discretization
  \(I_{B} = \{0=\theta_0< \theta_1, \dots, < \theta_{B - 1} < 2\pi\}\) of the
  interval \([0, 2\pi]\) (\Cref{sec:curvilinear}).
\item Using the discretization \(I_{B}\), construct two curvilinear
  meshes \(V^-_{B, d}\) and \(V^+_{B, C_r}\) in the near-boundary
  regions $V^-$ (within \(\overline{\Omega}\)) and $V^+$ (outside
  \(\Omega\)), respectively (equations~\eqref{eq:vnbj_m}
  and~\eqref{eq:vnbj_p}). Note that the discrete boundary points
  \(\bsym{q}(\theta_j)\) with \(\theta_j\in I_B\) are common to both
  \(V^-_{B,d}\) and \(V^+_{B,C_r}\). 
\item Determine the set \(H\cap V^+\) of Cartesian grid points and
  construct the proximity map
  \(\mathcal{P}: H\cap V^+ \to V^+_{B,C_r}\)
  (\Cref{sec:Proximity_map}).
\item
  For all \(Q\in H\cap V^+\) obtain the the parameter value
  \(\theta^Q\in [0,1]\) of the foot \(\mathcal{F}(Q)\) of the normal through
  \(Q\) (\Cref{sec:transf-fc-val-2-cartg}
  and the left image in \Cref{fig:interpolation_to_cartesian_grid}).  
\item For each normal grid line (inherent in \(V^-_{B,d}\) and
  \(V^+_{B,C_r}\)) given by the discretization \(I_B\), compute the
  blending-to-zero function values along that normal
  (\Cref{sec:interp-scheme-bound}).
\item For all the points \(Q\in H\cap V^+\), obtain the continuation
  function value at \(Q\) by local 2D interpolation
  (\Cref{sec:transf-fc-val-2-cartg} and the right image in
  \Cref{fig:interpolation_to_cartesian_grid}).
\item Apply the two-dimensional FFT once to the continuation
  function values to obtain the desired Fourier series in
  (\ref{eq:fc_2d_srs}). 
\end{enumerate}

\subsection{2D-FC approximation: Numerical Results}\label{sec:numer-results-fc2d}
\begin{figure}[h]
  \centering \captionsetup[captionsetup]{labelformat=empty}
  \begin{subfigure}[b]{0.35\textwidth}\centering
    \includegraphics[width=\linewidth]{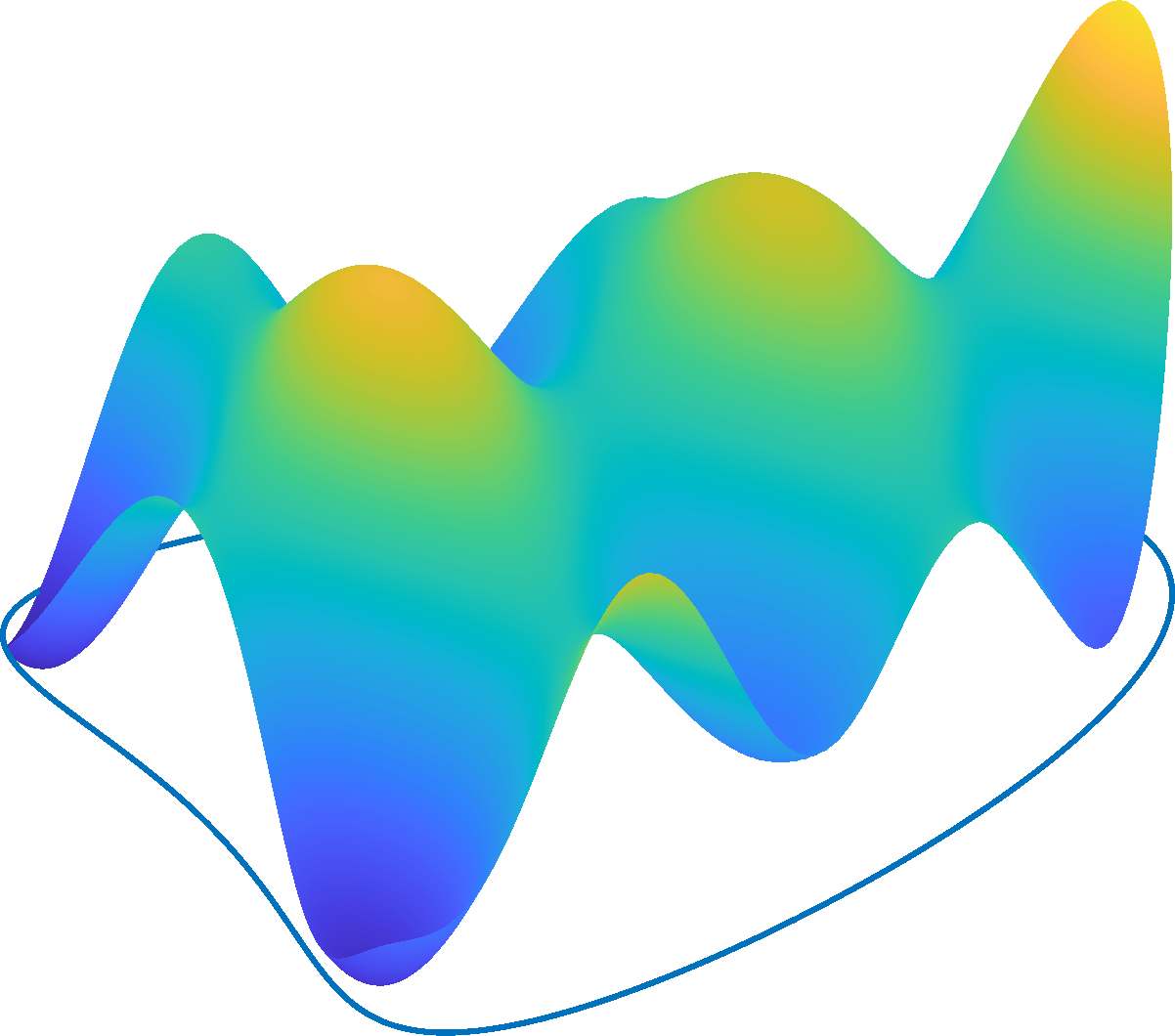}
  \end{subfigure}
  \hspace{1cm}
  \begin{subfigure}[b]{0.55\textwidth}   
    \includegraphics[width=\linewidth]{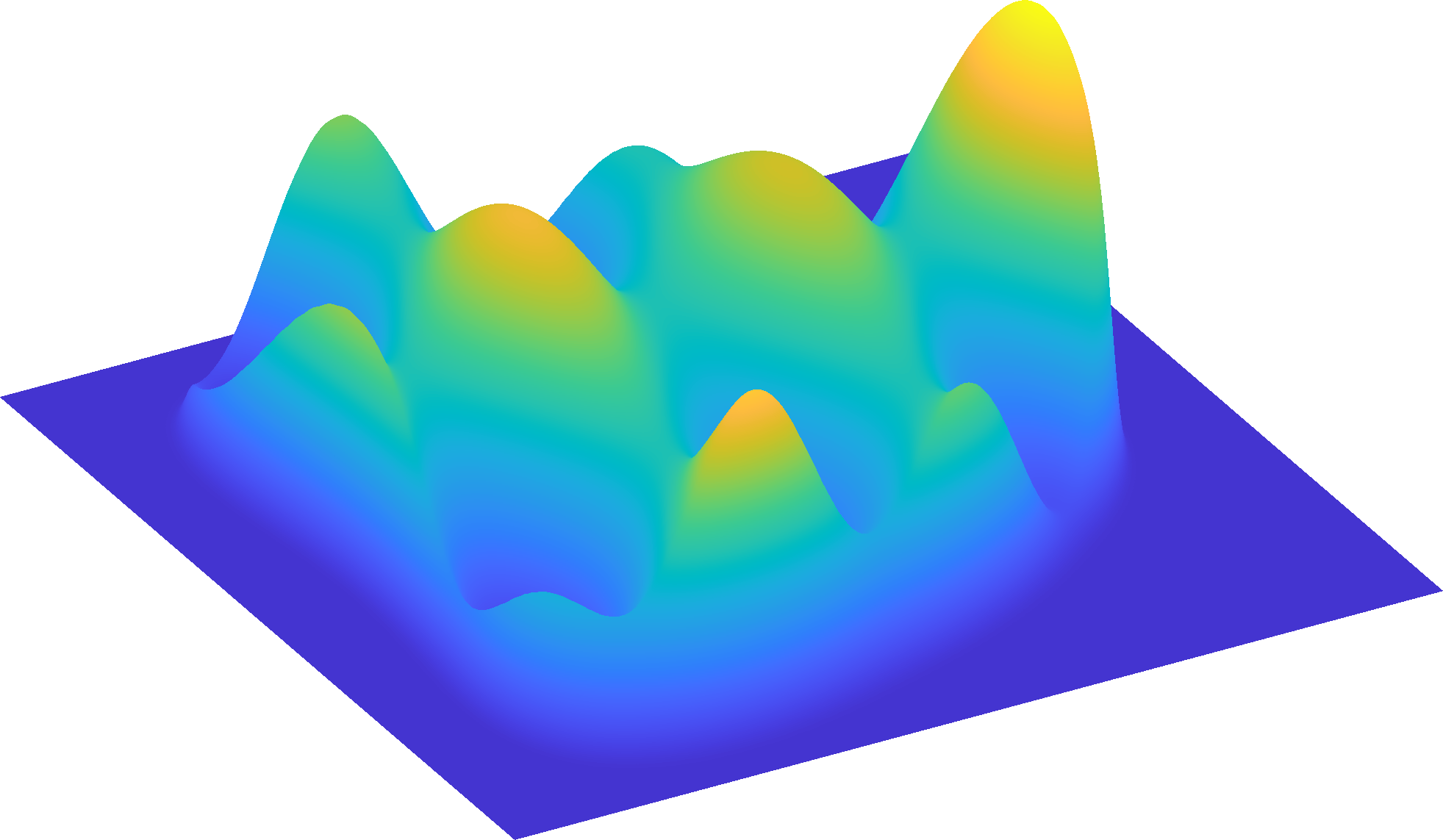}
  \end{subfigure}
  \caption{Demonstration of the 2D-FC procedure. Continuation of an
    oscillatory function defined on a kite shaped domain, as detailed
    in~\Cref{exm:fc_kite_pic}. The left and right images display the
    original and 2D-FC function values, respectively.  The blue curve
    in the left image indicates the boundary of the kite-shaped domain
    \(\Omega\). Note the narrowness of the region wherein the
    transition to zero takes place.}~\label{fig:demo_fc_kite}
\end{figure}

This section demonstrates the accuracy and efficiency of the proposed
2D-FC method. Use of the 2D-FC method requires selection of specific
values for each one of the following parameters (all of which are
introduced in \Cref{sec:fc_1d,sec:four-cont-meth}):
\begin{itemize}[noitemsep]
\item \(d\): number of points in the boundary section (\Cref{sec:background_fc_1d}).
\item \(C\): number of continuation points
  (\Cref{sec:background_fc_1d}).
\item \(Z\): number of zero matching points (\Cref{sec:background_fc_1d}).
\item \(n_{os}\): oversampling factor used in the oversampled
  matching procedure (\Cref{sec:background_fc_1d}).  
\item \(n_r\): refinement factor along the normal directions in
  \(V^+_{B,C_r}\) (\Cref{rmk:Cont_ref}).
\item \(R\): the smallest rectangle containing \(\Omega\cup V^+\)
  (\Cref{sec:fc-two-dimensions}).
\item \(N = N_x\times N_y\): number of points in the uniform spatial
  grid \(H\) (\Cref{sec:fc-two-dimensions}).
\item \(B\): number of points in the boundary
  discretization (\Cref{sec:fc-two-dimensions}). 
\item \(M-1\): interpolating polynomial degree
  (Sections~\ref{sec:interp-scheme-bound}
  and~\ref{sec:transf-fc-val-2-cartg}).
\end{itemize}
All the errors reported in this section were computed on a Cartesian
grid of step size $h/2$ within \(\Omega\).  In all of the numerical
examples considered in this article the parameter selections were made
in accordance with \Cref{rmk:param_sel}. The computer system used, in
turn, is described in \Cref{rmk:comp_sys}.
\begin{rembold}[Parameter selections]\label{rmk:param_sel}
  For a given step-size \(h\) in the two-dimensional Cartesian grid
  \(H\), the normal and the boundary step-sizes \(k_1\) and \(k_2\)
  (\Cref{sec:fc-two-dimensions}) were taken to coincide with $h$:
  $k_1=k_2=h$. The parameter values \(C = 27\), \(n_{os} = 20\),
  \(Z = 12\) and \(n_r = 6\) were used in the evaluation of the
  matrices \(\bsym{A}_r\) and \(\bsym{Q}\) (see
  \Cref{sec:background_fc_1d} and \Cref{rmk:Cont_ref}).  And, finally,
  with exception of the interpolation-degree experiments presented in
  \Cref{exm:fc_2d_conv_kite} (\Cref{tab:fc_2d_conv_kite_d4_d5}) and
  \Cref{exm:poisson_2} (\Cref{tab:poisson_conv_kite_d4}), the
  interpolating-polynomial degree \((M-1) = (d + 2)\) was used for the
  various matching-point numbers $d$ considered.
\end{rembold}
\begin{rembold}[Computer system]\label{rmk:comp_sys}
All of the numerical results reported in this paper
were run on a single core of a 3.40 GHz Intel Core i7-6700 processor
with 15.4Gb of 2133 MHz memory.
\end{rembold}
\begin{table}[h]
  \centering
  \begin{tabular}{c c c c c c c c}
    \toprule
    \(h\) & \(N^\Omega\) & \(B\) & \(T_{\mathcal{P}}\) & \(T_{\mathcal{F}}\) & \(T_V\) & \(T_H\)
    & Rel.~Err.~(\(\ell_\infty\))\\
    \midrule
    \(2\cdot 10^{-2}\) & 100 & 313 & \(2.1\cdot 10^{-3}\) & \(2.6\cdot 10^{-3}\) & \(1.7\cdot 10^{-3}\) & \(6.5\cdot 10^{-3}\) & \(4.7\cdot 10^{-4}\)\rdelim\}{3}{0mm}[\(d = 4\)] \\
    \(1\cdot 10^{-2}\) & 200 & 628 & \(4.2\cdot 10^{-3}\) & \(5.2\cdot 10^{-3}\) & \(3.8\cdot 10^{-3}\) & \(1.2\cdot 10^{-2}\) & \(2.7\cdot 10^{-5}\)\\
    \(5\cdot 10^{-3}\) & 400 & 1250 & \(1.0\cdot 10^{-3}\) & \(1.2\cdot 10^{-2}\) & \(1.0\cdot 10^{-4}\) & \(2.5\cdot 10^{-2}\) & \(1.8\cdot 10^{-6}\)\\
    \midrule
    \(2\cdot 10^{-2}\) & 100 & 313 & \(2.1\cdot 10^{-3}\) & \(2.9\cdot 10^{-3}\) & \(2.2\cdot 10^{-3}\) & \(9.0\cdot 10^{-3}\) & \(1.4\cdot 10^{-4}\)\rdelim\}{3}{0mm}[\(d = 5\)] \\
    \(1\cdot 10^{-2}\) & 200 & 628 & \(5.9\cdot 10^{-3}\) & \(6.9\cdot 10^{-3}\) & \(5.8\cdot 10^{-3}\) & \(1.9\cdot 10^{-2}\) & \(4.3\cdot 10^{-6}\)\\
    \(5\cdot 10^{-3}\) & 400 & 1250 & \(1.1\cdot 10^{-2}\) & \(1.3\cdot 10^{-2}\) & \(1.2\cdot 10^{-2}\) & \(3.1\cdot 10^{-2}\) & \(1.4\cdot 10^{-7}\)\\
    \midrule
    \(2\cdot 10^{-2}\) & 100 & 313 & \(2.2\cdot 10^{-3}\) & \(3.9\cdot 10^{-3}\) & \(4.2\cdot 10^{-3}\) & \(2.6\cdot 10^{-2}\) & \(4.1\cdot 10^{-5}\)\rdelim\}{2}{0mm}[\(d = 10\)] \\
    \(1\cdot 10^{-2}\) & 200 & 628 & \(4.3\cdot 10^{-3}\) & \(6.6\cdot 10^{-3}\) & \(8.8\cdot 10^{-3}\) & \(4.7\cdot 10^{-2}\) & \(1.6\cdot 10^{-7}\)\\
    \bottomrule
  \end{tabular}
  \caption{Times (in sec.) required by the various tasks in the 2D-FC
    algorithm in the setting of~\Cref{exm:fc_2d_disc_1}. The times
    reported were calculated as time-averages over \(10\) runs. The
    integer \(N^\Omega\) equals the number of spatial grid points used
    over the diameter of the disc.
  }\label{tab:time_for_fc_parts_d_4_5_10}
\end{table}
\begin{example}[Performance and efficiency of the 2D-FC
  method]\label{exm:fc_2d_disc_1}
  In our first example we consider a problem of FC approximation of
  the function $f: \Omega\to \mathbb{R}$ given by
  \begin{align}\label{eq:fc_example_1}
    f(x,y) &= -\sin(5\pi x)\sin(5\pi y)
  \end{align} 
  on the unit disc
  $\Omega =\{(x,y)\in \mathbb{R}^2: x^2+y^2 \leq 1 \}$. The left graph
  in \Cref{fig:fc_2d_err_time_4_5_10} displays the relative
  \(\ell_\infty\) and \(\ell_2\) errors, in \(\log\)-\(\log\) scale,
  obtained from 2D-FC approximations of the function $f$, for three
  different values of polynomial degree $d$ defined in
  Section~\ref{sec:interp-scheme-bound}, namely, \(d=4\), \(d = 5\)
  and \(d=10\) ---demonstrating the respective fourth, fifth and tenth
  orders of convergence expected. Higher rates of convergence, which
  are useful in some cases, can be achieved by using higher values of
  $d$, as demonstrated in the context of the Poisson solver in
  \Cref{sec:numerical-examples-poisson}.  The corresponding computing
  costs, including ``Geometrical Setup'' cost as well as the
  ``Cartesian Evaluation'' cost are presented in the right graph of
  \Cref{fig:fc_2d_err_time_4_5_10}. The Geometrical setup cost
  combines the setup time for the grids
  (Section~\ref{sec:curvilinear}) \({V}^-_{B,d}\), \({V}^+_{B,C_r}\)
  and \(H\); the time \(T_\mathcal{P}\) required for construction of
  the proximity map \(\mathcal{P}\) (\Cref{sec:Proximity_map}); and
  the time \(T_\mathcal{F}\) required for evaluation of the foot of
  the normal for all points in \(H\cap V^+\)
  (Section~\ref{sec:transf-fc-val-2-cartg}).  The Cartesian Evaluation
  time, in turn, equals the sum of the time \(T_V\) required for
  evaluation of the \(f^c\) values on the curvilinear grid \({V}^+\)
  and time \(T_H\) required for subsequent interpolation onto the
  Cartesian grid \(H\). \Cref{tab:time_for_fc_parts_d_4_5_10} reports
  additional details concerning computing times required by various
  tasks associated with the Geometrical Setup and Cartesian Evaluation
  for this example, including the times \(T_\mathcal{P}\),
  \(T_\mathcal{F}\), \(T_V\) and \(T_H\).  The Cartesian-interpolation
  time \(T_H\) dominates the overall Cartesian Evaluation step
  (cf.~\Cref{tab:time_for_fc_parts_d_4_5_10}).  In all of these cases
  we see that the computation time grows linearly with $1/h$. Also,
  the slope of the Cartesian evaluation cost depends on the degree
  \(d\) whereas the Geometrical setup cost, which remains similar in
  all the cases considered in this example, depends mainly on \(B\)
  and the refinement factor \(n_r\).
\end{example}
\begin{rembold}[Use of higher degree Gram polynomials]\label{rmk:higher_gram}
  Comparison of the various accuracy and timing values reported in
  Figure~\ref{fig:fc_2d_err_time_4_5_10} suggests that use of lower 2D
  FC degrees such as $d=4$ or $d=5$ may provide the highest efficiency
  for approximation accuracies up to single precision.
\end{rembold}
\begin{example}[Interpolation degree \((M-1)\) for a given
  2D-FC order \(d\)]\label{exm:fc_2d_conv_kite}
  Our next example concerns the approximation of the trigonometric
  function
  \begin{align}
    f &= -(x ^ 6 + y ^ 6)\sin(10\pi x)\sin(10\pi y), 
  \end{align}
  defined over the kite shaped domain contained within the curve given
  by \(x(\theta) = \cos(\theta) + 0.35\cos(2\theta) - 0.35\);
  \(y(\theta) = 0.7\sin(\theta)\) for \(0\leq \theta\leq 2\pi\).  A
  convergence study for this test case is presented
  in~\Cref{tab:fc_2d_conv_kite_d4_d5} for the values \(d = 4\) and
  \(d = 5\), and with \(M = d+1\), \(M = d+2\) and \(M = d+3\).
  Clearly, the selections $M=d+2$ and $M=d+3$ provide similar accuracy
  in most of the 2D-FC approximation cases considered. The value
  $M=d+3$, which yields somewhat better interpolation accuracy for
  larger step-sizes (e.g. \(h=10^{-2}\)), and which, as illustrated
  in~\Cref{tab:poisson_conv_kite_d4}, gives rise to some improvements
  for all step-sizes in the Poisson-problem applications considered in
  \Cref{sec:solut-poiss-equat}, was used in all of the numerical
  experiments presented in this article (except for the cases
  specifically designed to test the dependence of the accuracy on
  variations of the parameter $M$).
\end{example}
\begin{table}[h]
  \centering
  \begin{tabular}{c c c l c l c l c}
    \toprule
    & & & \multicolumn{2}{c}{\(M=d+1\)} & \multicolumn{2}{c}{\(M=d+2\)} & \multicolumn{2}{c}{\(M=d+3\)} \\
    \(h\) & \(N_x\) & \(N_y\) & Abs.~Err. & Order & Abs.~Err. & Order & Abs.~Err. & Order \\
    \midrule
    \(1\cdot 10^{-2}\) & \(297\) & \(197\) & \(1.8\cdot 10^{-3}\) & \:---\: & \(1.4\cdot 10^{-3}\) & \:---\: & \(9.2\cdot 10^{-4}\) & \:---\: \rdelim\}{5}{0mm}[\(d = 4\)]\\ 
    \(5\cdot 10^{-3}\) & \(537\) & \(337\) & \(2.0\cdot 10^{-4}\) & \(3.2\) &\(9.0\cdot 10^{-5}\) & \(3.9\) & \(3.1\cdot 10^{-5}\) & \(4.9\) \\
    \(2.5\cdot 10^{-3}\) & \(1017\) & \(617\) & \(1.0\cdot 10^{-5}\) & \(4.3\) &\(4.6\cdot 10^{-6}\) & \(4.3\) & \(2.3\cdot 10^{-6}\) & \(3.8\) \\
    \(1.25\cdot 10^{-3}\) & \(1977\) & \(1177\) &\(6.2\cdot 10^{-7}\) & \(4.0\) &\(1.4\cdot 10^{-7}\) & \(5.0\) & \(1.4\cdot 10^{-7}\) & \(4.0\) \\
    \(6.25\cdot 10^{-4}\) & \(3897\) & \(2297\) &\(3.4\cdot 10^{-8}\) & \(4.2\) &\(9.0\cdot 10^{-9}\) & \(4.0\) & \(9.0\cdot 10^{-9}\) & \(4.0\) \\
    \midrule
    \(1\cdot 10^{-2}\) & \(297\) & \(197\) & \(4.4\cdot 10^{-3}\) & \:---\: &\(2.3\cdot 10^{-3}\) & \:---\: & \(2.5\cdot 10^{-4}\) & \:---\:  \rdelim\}{5}{0mm}[\(d = 5\)] \\    
    \(5\cdot 10^{-3}\) & \(537\) & \(337\) & \(3.3\cdot 10^{-4}\) & \(3.7\) &\(4.2\cdot 10^{-5}\) & \(5.8\) & \(1.5\cdot 10^{-5}\) & \(4.1\) \\
    \(2.5\cdot 10^{-3}\) & \(1017\) & \(617\) & \(1.8\cdot 10^{-5}\) & \(4.2\) &\(4.3\cdot 10^{-7}\) & \(6.6\) & \(2.6\cdot 10^{-7}\) & \(5.8\) \\
    \(1.25\cdot 10^{-3}\) & \(1977\) & \(1177\) & \(4.1\cdot 10^{-7}\)& \(5.4\) &\(6.4\cdot 10^{-9}\) & \(6.1\) & \(4.1\cdot 10^{-9}\) & \(6.0\) \\
    \(6.25\cdot 10^{-4}\) & \(3897\) & \(2297\) & \(1.2\cdot 10^{-8}\)& \(5.0\) &\(1.3\cdot 10^{-10}\) & \(5.6\) & \(1.3\cdot 10^{-10}\) & \(5.0\) \\
    \bottomrule
  \end{tabular}
  \caption{Convergence table for the 2D-FC method in the setting
    of~\Cref{exm:fc_2d_conv_kite}. For both \(d=4\) and \(d=5\), we
    observe the expected fourth and fifth orders of convergence,
    respectively, for all the three choices of \(M\), namely,
    \(M=d+1\), \(M=d+2\) and \(M=d+3\). The value \(M=d+3\) leads to
    somewhat improved accuracy. }~\label{tab:fc_2d_conv_kite_d4_d5}
\end{table}
\begin{example}[Graphical illustration of the 2D-FC
  method]\label{exm:fc_kite_pic}
  \Cref{fig:demo_fc_kite} demonstrates the 2D-FC extension method for
  the function defined by
  \[
    f(x, y) = 4+(1+x^2+y^2)(\sin(2.5\pi x-0.5)+\cos(2\pi y-0.5)),
  \]
  over the kite shaped domain considered in
  \Cref{exm:fc_2d_conv_kite}. Both the original function and its
  extension are presented in \Cref{fig:demo_fc_kite}.
\end{example}

\section{Applications of the 2D-FC method}\label{sec:applications-2d-fc}
The 2D-FC method introduced in Section~\ref{sec:fc-two-dimensions} can
be used to facilitate spectral treatment in cases for which iterated
use of the 1D Fourier expansions does not suffice, but for which use
of full 2D Fourier expansions is
beneficial. Sections~\ref{sec:solut-poiss-equat} and~\ref{sec:fc_ff}
briefly describe two such cases, one concerning solution of the
Poisson problem via Fourier expansions, and, the other, the solution
of the wave equation via a novel Fourier Forwarding (FF) technique. A
few numerical examples are presented illustrating the character of the
resulting numerical solvers.

\subsection{Application Example I:~Poisson Problem\label{sec:solut-poiss-equat}}
In this section we present a 2D-FC based method for the solution of
the 2D Dirichlet Poisson problem
\begin{equation}\label{eq:def_poisson}
  \begin{cases}
    \Delta u(x, y) = f(x,y), &(x, y)\in\Omega\\
    u(x, y) = g(x, y), &(x, y)\in \Gamma;
  \end{cases}
\end{equation}
the corresponding problem under Neumann or Robin boundary conditions
can be treated similarly. Only a brief sketch is presented here, as an
illustration of the 2D-FC approach.  A complete description of the
method, including a detailed geometrical treatment needed for reliable
evaluation of the solution near boundaries, is presented
in~\cite{fc_based_poisson_solver}.

The proposed Poisson solver obtains the solution $u$ within the
prescribed tolerance as a sum
\begin{equation}\label{eq:total_soln}
  u = u_p + v
\end{equation}
of a ``particular solution'' $u_p$, produced by means of the 2D-FC
method, which satisfies the Poisson equation $\Delta u = f$ (but which
generically does not satisfy the boundary conditions), and a solution
$v$ of the ``homogeneous problem'', produced by means of a boundary
integral equation, which satisfies the Dirichlet boundary value
problem for Laplace's equation
\begin{equation}
  \label{eq:def_laplace}
  \begin{cases}
    \Delta v(x, y) = 0, &(x, y)\in\Omega\\
    v(x, y) = g_{\hom}(x, y), &(x, y)\in \Gamma, 
  \end{cases}
\end{equation}
where
\begin{equation}\label{eq:laplace_bc}
  g_{\hom}(x, y) = g(x, y) - u_p(x, y)|_{\Gamma}.
\end{equation}

A particular solution $u_p$ for the problem~\eqref{eq:def_poisson} can
easily be obtained from the 2D-FC expansion $f^c(x,y)$ of the
right-hand function $f(x,y)$ (equation~\eqref{eq:fc_2d_srs})---in view
of the diagonal character of the Laplace operator in Fourier space.
We thus obtain
\begin{align}\label{eq:poisson_ps}
  u_p(x, y) &= - \hat{f}^c_{0, 0} (x ^ 2 + y ^ 2) / 4 +
              \sum_{\ell = -N_x / 2 + 1} ^ {N_x / 2} \sum_{m = -N_y / 2 + 1} ^ {N_y / 2}
              b_{\ell,m}e ^ {2\pi i\left(\frac{\ell x}{L_x} + \frac{m y}{L_y}\right)}
\end{align}
where
\begin{align}\label{eq:poisson_ps_coeffs}
  b_{\ell,m}
  &= \begin{cases}
    \hfil 0, & \text{ if } (\ell, m) = (0, 0)\\
    \frac{ -  \hat{f}^c_{\ell, m}}{(2\pi \ell/ L_x)^2 + (2\pi m / L_y)^2},
    & \textrm{ if } (\ell, m)\neq (0, 0),
  \end{cases}
\end{align}
where one of a variety of possible selections was made for the
constant Laplacian term.  In view the asymptotically small factors
that relate the Fourier coefficients $b_{\ell,m}$ to the original FC
coefficients $\hat{f}^c_{\ell, m}$ it can be
shown~\cite{fc_based_poisson_solver} that, as illustrated
in~\Cref{sec:numerical-examples-poisson}, the rate of convergence in
the overall numerical 2D-FC based solution $u$ is of
\(\mathcal{O}(h^{d+2})\) if a 2D-FC algorithm of \(\mathcal{O}(h^d)\)
is utilized to compute the particular solution \(u_p\) (provided a
sufficiently accurate method is subsequently used for evaluation of
the homogeneous solution).

Values $u_p(\bsym{r})$ of the particular solution at points $\bsym{r}$
on the boundary $\Gamma$ are required as an input
(via~\eqref{eq:laplace_bc}) in the boundary-value
problem~\eqref{eq:def_laplace} for the Laplace solution $v$. It is
therefore necessary to utilize an efficient method for evaluation of
$u_p$ at points $\bsym{r}$ that are not part of the Cartesian mesh
$H$.  The straightforward procedure based on direct addition of all
terms in~\eqref{eq:poisson_ps} for each discretization point on
$\Gamma$ does not match the optimal $\mathcal{O}(N\log(N))$ cost
asymptotics enjoyed by all the other elements of the algorithm and is
therefore avoided.  Instead, the proposed algorithm first obtains the
values of the right hand side of~\eqref{eq:poisson_ps} for all
$\bsym{r}\in H$ via a direct application of the FFT algorithm, and,
then, using these values, it produces the values for
$\bsym{r}\in \Gamma$ via iterated one-dimensional interpolation, as
described in~\cite[Sec. 3.6.1]{numer-recipes}. In order to match the
overall order \((d+2)\) accuracy of the overall Poisson solution,
one-dimensional polynomial interpolants of degree \((M_P-1) \geq (d+2)\)
(cf.~\Cref{exm:poisson_2}) are used in this context for both the $x$
and $y$ interpolation directions.

The numerical solution of the Laplace equation
in~\eqref{eq:def_laplace}, in turn, can be obtained rapidly and
efficiently on the basis of the boundary integral method (see
e.g.~\cite{book-kress-lie-2014}). Relying on the boundary
parametrization~\eqref{eq:param}, the proposed algorithm incorporates
an integral equation with a smooth kernel together with the simple and
effective Nystr\"om algorithm presented
in~\cite[Sec. 12.2]{book-kress-lie-2014}. Based on trapezoidal-rule
quadrature, this algorithm results in highly accurate solutions: in
view of the periodicity and smoothness of the solution and the kernel,
the approach yields super-algebraically small errors provided the
boundary and right-hand side $g_{\hom}$ are both smooth. The
associated linear system is solved by means of the iterative linear
algebra solver \emph{GMRES}~\cite{gmres}. Note that the integrand
exhibits a near singular behavior for evaluation points that are near
the boundary $\Gamma$ but that are not on $\Gamma$. In order to
address this difficulty, the proposed method uses a scheme (some
elements of which were introduced
in~\cite{AKHMETGALIYEV20151,Bruno_Delourme_2014}) which, based on
local mesh refinement and subsequent interpolation using polynomial of
degree \((M_P-1)\), successfully resolves this difficulty.  A detailed
description of this and other aspects concerning the 2D-FC based
Poisson solver is presented in the forthcoming
contribution~\cite{fc_based_poisson_solver}.

Once the particular and homogeneous solutions $u_p$ and $v$ have been
obtained, the solution $u$ of the Poisson problem is given
by~\eqref{eq:total_soln}. The numerical convergence rate of the
solution produced by the algorithm is mainly determined by the order
$d$ of the 2D-FC algorithm used. In all, the method is fast and highly
accurate; a few illustrations, including accuracy and timing
comparisons with leading solvers, are presented in the following
section.
\subsubsection{Numerical Illustrations for the Poisson
  Problem\label{sec:numerical-examples-poisson}}
The numerical illustrations presented in this section utilize the
2D-FC parameter selections presented in~\Cref{rmk:param_sel}, with
various choices of the order parameter \(d\). In addition, the size of
the uniform boundary discretization used by the trapezoidal-rule based
Nystr\"{o}m method is taken, for simplicity, to equal \(N_x\)---but,
of course, in view of the super-algebraic convergence of the
trapezoidal-rule quadrature, a smaller discretization size could have
been used without sacrificing accuracy. The $\ell_2$ and $\ell_\infty$
errors reported in this section were computed over the Cartesian grid
\(H\cap \Omega\) unless indicated otherwise.  The first Poisson-solver
example concerns a simple problem previously considered
in~\cite{FRYKLUND201857}.

\begin{figure}[h]
  \centering
  \begin{subfigure}[b]{0.45\textwidth}  
    \includegraphics[scale=0.64]{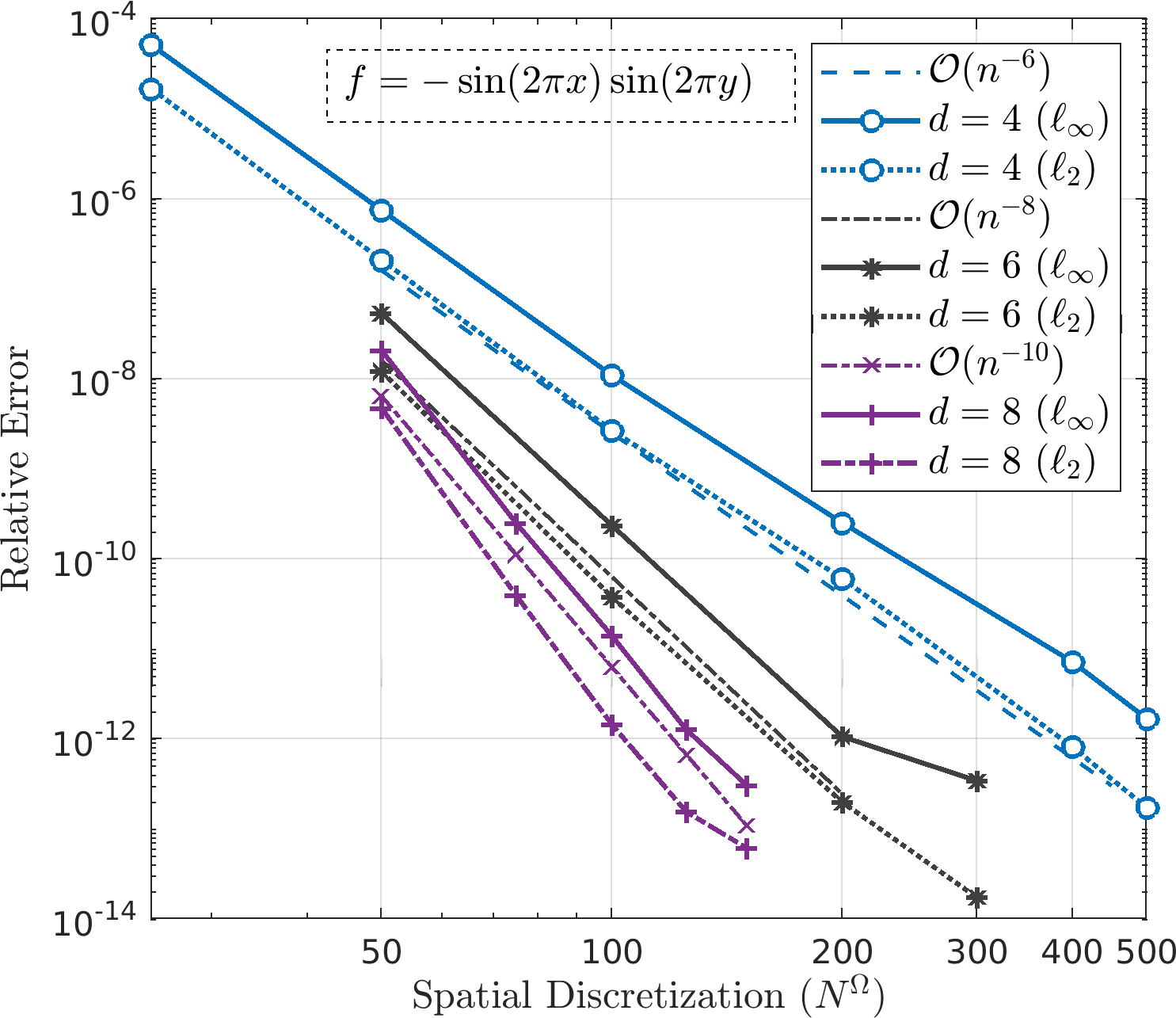}
  \end{subfigure}
  \begin{subfigure}[b]{0.45\textwidth}
    \includegraphics[scale=0.64]{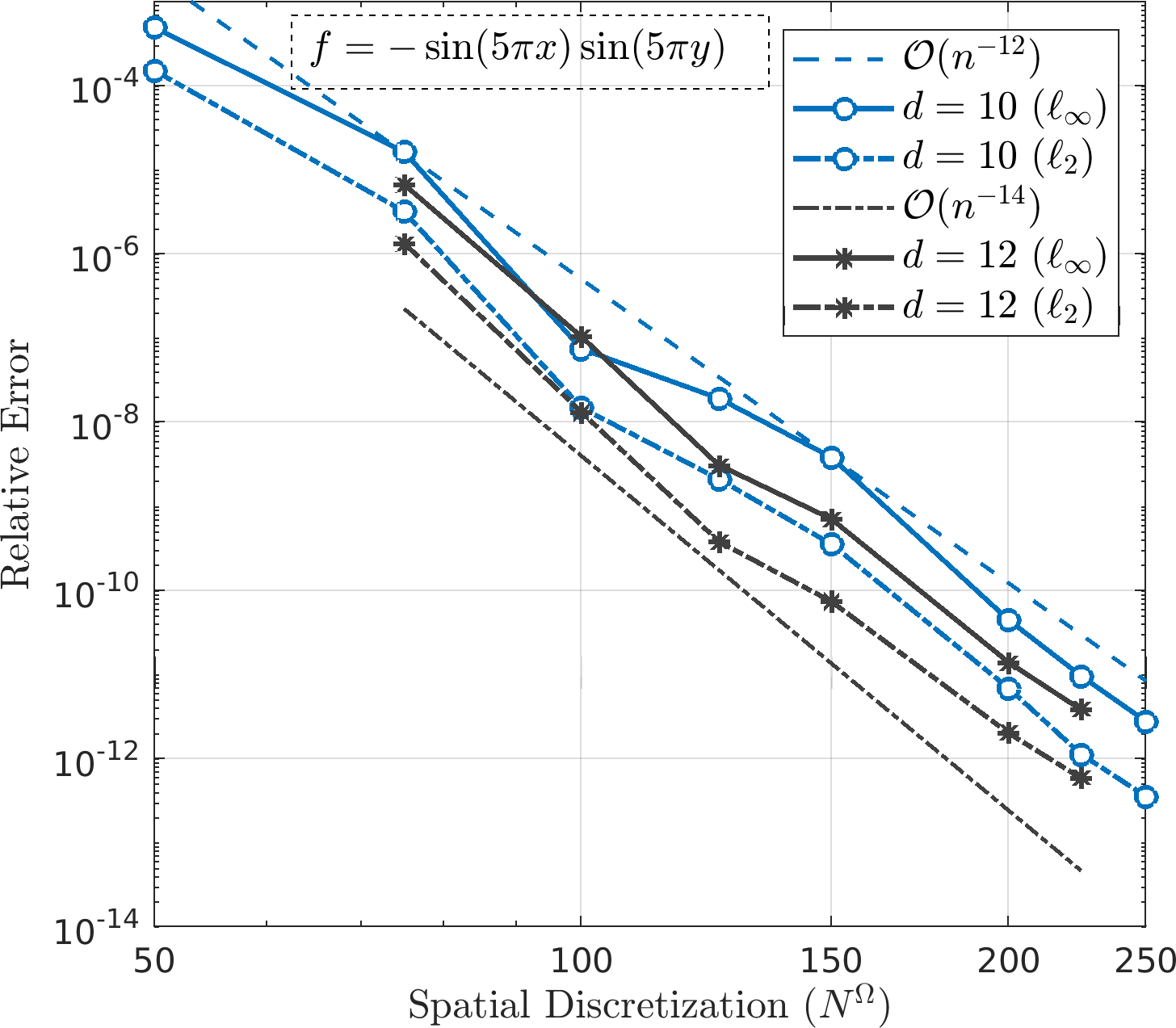}
  \end{subfigure}
  \caption{Numerical solution errors, in \(\log\)-\(\log\) scale,
    resulting from use of the 2D-FC based Poisson solver over the unit
    disc (\Cref{exm:poisson_1}). Left graph: solution errors for 2D-FC
    methods of orders \(d=4\), \(d=6\) and \(d=8\) for the function
    \(f = -\sin(2\pi x) \sin(2\pi y)\). Right graph: solution errors
    with \(d=10\) and \(d=12\) for the function
    \(f = -\sin(5\pi x)\sin(5\pi y)\). The parameter \(N^\Omega\)
    denotes the number of spatial grid points used over the diameter
    of the disc. A total of \(N_x=N_y=N^\Omega+2C\) points (\(C=27\))
    were used over each dimension of the periodic square $R$. As
    discussed in the text, a convergence rate of order $(d+2)$ results
    from use of a 2D-FC approximation of order
    $d$.}~\label{fig:poisson_errors}
\end{figure}

\begin{example}[High-order 2D-FC based Poisson solution]\label{exm:poisson_1}
  We consider the Poisson problem~\ref{eq:def_poisson} in the domain
  $\Omega = \{(x,y)\in\mathbb{R}^2 : x^2 + y^2 \leq 1\}$ with
  \(f = -\sin(2\pi x) \sin(2\pi y)\). The left portion of
  \Cref{fig:poisson_errors} presents the numerical errors in the
  solutions produced by the 2D-FC based Poisson solvers for \(d = 4\),
  \(d=6\) and \(d=8\) for \(f = -\sin(2\pi x) \sin(2\pi y)\). The
  observed rates of convergence for all the three cases match the
  expected increased rates of convergence, as discussed in
  \Cref{sec:solut-poiss-equat}, that is, rates convergence of orders
  \(6\), \(8\) and \(10\), respectively.  This problem was also
  considered in~\cite{FRYKLUND201857}.  Comparison of the results
  presented in~\cite[Fig.~8]{FRYKLUND201857} and those on the left
  graph in~\Cref{fig:poisson_errors} suggests the 2D-FC based Poisson
  solver performs favorably for high accuracies. For instance, a
  number \(N^\Omega=100\) of spatial grid points over the diameter of
  \(\Omega\), that is to say, \(N_x=154\) grid points over one length
  of the rectangular computational domain, provides, as shown in
  the~\Cref{fig:poisson_errors}, an \(\ell_2\) error
  \(1.4\cdot 10^{-12}\) whereas, in~\cite[Fig.~8]{FRYKLUND201857}, a
  similar discretization provides \(\ell_2\) errors close to
  \(10^{-9}\). The \(10^{-12}\) error is achieved in that reference at
  a number of approximately \(275\) points in spatial discretization
  points in each spatial direction.  For a different test case in this
  problem setting we now take \(f = -\sin(5\pi x) \sin(5\pi y)\) (a
  function that was also used for the convergence study of the 2D-FC
  algorithm as presented in \Cref{exm:fc_2d_disc_1}), and we report,
  on the right graph in \Cref{fig:poisson_errors}, the numerical
  errors in the solution produced by the solvers for higher values of
  \(d\), namely, \(d=10\) and \(d=12\). Once again the expected
  convergence rates (in this case, of orders \(12\) and \(14\),
  respectively) are observed in practice.
\end{example}
\begin{table}[h]
  \centering
  \begin{tabular}{c c c l c l c l c}
    \toprule
    & & & \multicolumn{2}{c}{\(M=M_P=d+1\)} & \multicolumn{2}{c}{\(M=M_P=d+2\)} & \multicolumn{2}{c}{\(M=M_P=d+3\)} \\
    \(h\) & \(N_x\) & \(N_y\) & Rel.~Err.(\(\ell_2\)) & Order & Rel.~Err.(\(\ell_2\)) & Order & Rel.~Err.(\(\ell_2\)) & Order \\
    \midrule
    \(4\cdot 10^{-2}\) & \(117\) & \(93\) & \(6.0\cdot 10^{-6}\) & \:---\: & \(7.5\cdot 10^{-7}\) & \:---\: & \(4.1\cdot 10^{-7}\) & \:---\:\\
    \(2\cdot 10^{-2}\) & \(177\) & \(127\) & \(1.7\cdot 10^{-7}\) & \(5.1\) &\(1.0\cdot 10^{-8}\) & \(6.2\) & \(6.3\cdot 10^{-9}\) & \(6.0\) \\
    \(1\cdot 10^{-2}\) & \(297\) & \(197\) & \(5.2\cdot 10^{-9}\) & \(5.0\) &\(1.5\cdot 10^{-10}\) & \(6.1\) & \(8.7\cdot 10^{-11}\) & \(6.2\) \\
    \(5\cdot 10^{-3}\) & \(537\) & \(337\) & \(1.6\cdot 10^{-10}\) & \(5.1\) &\(2.6\cdot 10^{-12}\) & \(5.9\) & \(1.2\cdot 10^{-12}\) & \(6.2\) \\
    \(4\cdot 10^{-3}\) & \(657\) & \(407\) & \(5.2\cdot 10^{-11}\) & \(4.9\) &\(6.8\cdot 10^{-13}\) & \(6.0\) & \(3.0\cdot 10^{-13}\) & \(6.2\) \\
    \bottomrule
  \end{tabular}
  \caption{Convergence of the 2D-FC based solution of the Poisson
    problem described in~\Cref{exm:poisson_2}.
  }~\label{tab:poisson_conv_kite_d4}
\end{table}
\begin{example}[Poisson solution interpolation degree \((M_P-1)\)]\label{exm:poisson_2}
  Once again we consider the problem~\eqref{eq:def_poisson} over a
  kite shaped domain as in~\Cref{exm:fc_2d_conv_kite} with
  \(f = -\sin(2\pi x) \sin(2\pi y)\). The errors in the solutions
  produced by the 2D-FC based Poisson solver and the corresponding
  convergence rates for \(d=4\) and three different values of \(M_P(=M)\),
  namely, \(M_P=d+1\), \(M_P=d+2\) and \(M_P=d+3\), are presented
  in~\Cref{tab:poisson_conv_kite_d4}. The observed rates of
  convergence for \(M_P=d+2\) and \(M_P=d+3\) show the increased
  \(d+2=6\)-th order convergence rate whereas the selection
  \(M_P=d+1=5\) shows a fifth order convergence as the overall error
  in the Poisson solution is dominated by the error associated with
  the order-five interpolation process. The value \(M_P=d+3\) is
  utilized for all Poisson-problem numerical results presented in this
  paper.
\end{example}
\begin{example}[Highly oscillatory Poisson problem]\label{exm:poisson_3}
  Here we consider the problem~\eqref{eq:def_poisson} over the kite
  shaped domain considered in~\Cref{exm:fc_2d_conv_kite} with the
  highly oscillatory right hand side
  \(f = -\sin(40\pi x) \sin(40\pi y)\). In this example, where we have
  used \(d=10\) for the 2D-FC particular-solution algorithm, the
  overall convergence rate, as reported
  in~\Cref{tab:conv_poisson_2_d10}, is close to the expected
  convergence rate of order $(d+2)$. In order to avoid near-singular
  integration problems which arise, under the fine discretizations
  considered in this example, as the numerical solution is evaluated
  at points very near the boundary $\Gamma$, here we report the error
  at all the Cartesian points within the computational domain that are
  at a greater distance from \(\Gamma\) than \(0.2\). Near boundary
  evaluation algorithms capable of evaluating the solution at points
  arbitrarily close to $\Gamma$ are presented
  in~\cite{fc_based_poisson_solver}.
\end{example}
\begin{table}[h]\label{tab:conv_poisson_3}
  \centering
  \begin{tabular}{c c c l c l c}
    \toprule
    \(h\) & \(N_x\) & \(N_y\) & Rel.~Err.~(\(\ell_\infty\)) & Order & Rel.~Err.~(\(\ell_2\)) & Order \\
    \midrule 
    \(5\cdot 10^{-3}\) &\(537\) & \(337\) & \(9.3\cdot 10^{-4}\) & --- &\(3.4\cdot 10^{-4}\) & --- \\
    \(2.5\cdot 10^{-3}\)&\(1017\) & \(617\)  & \(6.0\cdot 10^{-8}\) & \(13.9\) &\(1.7\cdot 10^{-8}\) & \(14.3\)\\
    \(1.25\cdot 10^{-3}\)&\(1977\) & \(1177\) & \(5.8\cdot 10^{-12}\) & \(13.3\) & \(1.6\cdot 10^{-12}\) & \(13.4\)\\
    \bottomrule
  \end{tabular}
  \caption{Interior errors in the numerical solutions produced by the
    2D-FC based Poisson solver of order \(d=10\), in the setting
    of~\Cref{exm:poisson_3}. (Errors are computed over points lying at
    a distance \(0.2\) from the domain boundary;
    see~\Cref{exm:poisson_3}.) A scaling-error even better than the
    expected order $(d+2) = 12$ was observed in this case.
  }\label{tab:conv_poisson_2_d10}
\end{table}

\subsection{Application Example II:~Fourier Forwarding (FF) method for Wave propagation problems\label{sec:fc_ff}}
This section presents the 2D-FC based Fourier-Forwarding method (FF)
for the solution of the wave equation and other constant coefficient
hyperbolic problems. Only a brief sketch of the FF approach, for
demonstration purposes, is presented here; a detailed account of this
methodology, including a variety of techniques designed for treatment
of boundary conditions, and with application to linear hyperbolic
systems, including, e.g., treatment of propagation in anisotropic
media, is presented in the forthcoming
contribution~\cite{fc_based_wave_solver_ff}.  In what follows
attention is restricted to the initial boundary value problem for the
wave equation in 2D, namely
\begin{equation}\label{eq:wave_2d}
  u_{tt} = c^2 (u_{xx} + u_{yy}), \quad\textrm{ for }
  (x, y,t ) \in \Omega\times\mathbb{R}_+
\end{equation}
with initial conditions \(u(x,y,t)|_{t=0} = f(x, y)\) and
\(u_t(x,y,t)|_{t=0} = g(x, y)\) for $(x, y)\in\Omega$, and with
appropriate boundary conditions on $\Gamma$.

In order to obtain a solution of this problem, the FF method
1)~Exploits the fact (also used in~\cite{VAY2013260,haber_psatd_1973}
in the context of bi-periodic problems) that the solution of the wave
equation in all of $\mathbb{R}^2$ with the initial data
$u(\bsym{r},0) = a e^{i \bsym{\kappa}\cdot\bsym{r}}$ and
$u_t(\bsym{r},0) = b e^{i \bsym{\kappa}\cdot\bsym{r}}$ is itself given
in closed form as a combination of two time-domain plane waves; and,
2)~Constructs auxiliary solutions $U_{FC}(x,y,t;T_j)$ of the form
\begin{equation}\label{eq:aux_sol}
  U_{FC}(x,y,t;T_j)
  =  \sum_{\ell = -N_x / 2 + 1}^{N_x / 2}\sum_{m = -N_y / 2 + 1}^{N_y / 2}
  a_{\ell m}(t;T_j) e ^ {2\pi i\left(\frac{\ell x}{L_x} + \frac{m y}{L_y}\right)},
\end{equation}
of~\cref{eq:wave_2d} utilizing the 2D-FC algorithm on certain initial
solution values at various times: $T_j = j \Delta T$
($j = 0,\dots, n$) with a ``large'' time step value $\Delta T$ (that
is selected so as to optimize the overall computational cost of the FF
algorithm), and for an arbitrary user-prescribed positive integer
$n$. In view of the limited domain of dependence of solutions of the
wave equation~\cite{John_1982}, the auxiliary solution
$U_{FC}(x,y,t;T_j)$, for \(T_j\leq t \leq T_{j+1}\), provides a valid
numerical approximation of the solution \(u(x,y,t)\) over a certain
subset
\(\Omega_{\Delta T}=\{\bsym{r} = (x,y)\in\Omega :
\mbox{dist}(\bsym{r},\Gamma) \geq c \Delta T\}\), away from the
boundary \(\Gamma\), of the domain \(\Omega\). To compute the solution
$U_B(x,y,t;T_j)$ on the region \(\Omega\setminus\Omega_{\Delta T}\)
adjacent to the physical boundary, the FF method uses a classical
time-stepping scheme, with spatial derivatives obtained by means of
the 1D-FC method, and using a (typically much smaller) time step
\(\Delta t\), which should be small enough so to ensure stability (as
dictated by the CFL condition) and to yield an accuracy level
consistent with that inherent in the 2D-FC approximation used. In
what follows we discuss the evaluation of the solution
\begin{equation}\label{eq:ff_soln}
  u(x,y,t) =
  \begin{cases}
    U_{FC}(x,y,t;T_0) &\textrm{ for } (x,y,t)\in \Omega_{\Delta T}\times [T_0,T_1],\\
    U_B(x,y,t;T_0) &\textrm{ for } (x,y,t)\in
    \Omega\setminus\Omega_{\Delta T}\times [T_0,T_1]
  \end{cases}
\end{equation}
for the time interval \(T_0=0\leq t \leq T_1=\Delta T\); a similar
procedure can be used to evaluate, inductively, the solution at all
other time intervals \(T_{j}\leq t \leq T_{j+1}\). For notational
simplicity the argument \(T_0\) is suppressed in what follows.

In order to obtain the auxiliary solution $U_{FC}(x,y,t)$ the method
utilizes the Cartesian grid $H$ (\cref{eq:cart_grid}) together with
the 2D-FC expansions
\begin{equation}\label{eq:wave_2d_ff_fc}
  \begin{cases}
    F(x,y) = \sum_{\ell = -N_x / 2 + 1}^{N_x / 2}\sum_{m = -N_y / 2 +
      1}^{N_y / 2}
    \hat{f}^c_{\ell m} e ^ {2\pi i\left(\frac{\ell x}{L_x} + \frac{m y}{L_y}\right)},\\
    G(x,y) = \sum_{\ell = -N_x / 2 + 1}^{N_x / 2}\sum_{m = -N_y / 2 +
      1}^{N_y / 2} \hat{g}^c_{\ell m} e ^ {2\pi i\left(\frac{\ell
          x}{L_x} + \frac{m y}{L_y}\right)},
  \end{cases}
\end{equation}
(cf.~\cref{eq:fc_2d_srs}) of the initial data \(f(x,y)\) and
\(g(x,y)\). (Note that while $F$ and $G$ are obtained from the given
initial conditions in the present case $j=0$, they are produced from
the numerical values of the solution $u(x,y,t)|_{t=T_j}$ and its time
derivative $u_t(x,y,t)|_{t=T_j}$ in the case $j>0$.)  Clearly,
provided the functions $a_{\ell m}(t)$ in~\eqref{eq:aux_sol} satisfy
the equations
\begin{align}\label{eq:wave_2d_ff_modal}
  \begin{cases}
    a ^ {\prime\prime}_{\ell m}(t) + \alpha_{\ell m} a_{\ell m}(t) = 0,\\
    a_{\ell m}(t)|_{t=0} = \hat{f}^c_{\ell m},\\
    \frac{\partial a_{\ell m}(t)}{\partial t}|_{t=0} =\hat{g}^c_{\ell m}
  \end{cases}
\end{align}
for $-N_x/2+1\leq \ell\leq N_x/2$ and $-N_y/2+1\leq m\leq N_y/2$
(where \(\alpha_{\ell m} = (2\pi c)^2[(\ell / L_x)^2 + (m/L_y)^2]\)),
the function $U_\mathrm{FC}(x,y,t)$ satisfies~\eqref{eq:wave_2d} as
well as the initial conditions $U_\mathrm{FC}(x,y,t)|_{t=0} = F(x,y)$
and
$\frac{\partial U_\mathrm{FC}(x,y,t)}{ \partial t}|_{t=0} = G(x,y)$
for \((x,y)\in \mathbb{R}^2\). Substituting the explicit solutions
\begin{align}\label{eq:wave_2d_ff}
  a_{\ell m}(t)
  &= \begin{cases}
    \hat{f}^c_{\ell m} + t \hat{g}^c_{\ell m}, & \textrm{ for } (\ell, m) = (0, 0)\\
    \hat{f}^c_{\ell m}\cos(\alpha_{\ell m}t)
    + \dfrac{\hat{g}^c_{\ell m}}{\alpha_{\ell m}}\sin(\alpha_{\ell m} t),
    & \textrm{ for } (\ell, m)\neq (0, 0). 
  \end{cases}
\end{align}
of the ODE~\eqref{eq:wave_2d_ff_modal} into~\eqref{eq:aux_sol} the
solution \(U_\mathrm{FC}(x,y,t)\) for $(x,y)\in\mathbb{R}^2$ is
obtained. An application of the two-dimensional spatial inverse FFT
over the Cartesian grid $H$ to the coefficients \(a_{\ell,m}(T_1)\)
then yields, per the discussion above concerning domains of
dependence, a numerical approximation of the solution \(u(x,y,t)\) for
\(T_0\leq t \leq T_1\) and for all
\((x,y)\in H\cap \Omega_{\Delta T}\). Note that the accuracy of the
auxiliary solution in its domain of validity \(\Omega_{\Delta T}\) at
time \(t\) (\(T_0\leq t \leq T_1\)), is only limited by the accuracy of
the underlying 2D-FC approximation of the functions $f$ and $g$ by the
FC expansions $F$ and $G$ throughout $\Omega$, respectively.

To obtain the near boundary solution \(U_\mathrm{B}(x,y,t)\), on the
other hand, the method uses a classical explicit time stepping scheme
in a certain open set
$\Omega_B\supset\Omega\setminus \Omega_{\Delta T}$ adjacent to the
boundary. In the proposed near-boundary evaluation algorithm (which is
based on use of the small time step $\Delta t$ on a Cartesian grid on
$\Omega_B$) various time stepping schemes, including the
Adams-Bashforth~\cite{fdm_leveque} and the Taylor
series~\cite{appelo-peterson-2012} methods, can be utilized; the
required spatial derivatives, in turn, are computed, with spectral
accuracy and without dispersion, by means of the 1D-FC method
(\Cref{sec:background_fc_1d}; cf.
also~\cite{amlani-bruno-fc-spectral-2016-307,BRUNO20102009}). In the
present $j=0$ case the initial values for \(U_B\) are obtained from
the initial conditions \(f(x,y)\) and \(g(x,y)\); for subsequent time
intervals the solution process for \(U_B\) is simply continued forward
in time: no additional initial values are needed for $U_\mathrm{B}$ at
the start of the time intervals $[T_j,T_{j+1}]$ for $j>0$. Boundary
conditions for $U_\mathrm{B}$ must be enforced at all boundary points
in $H \setminus\Omega_{\Delta T}$---including those near $\Gamma$ and
those near $\partial\Omega_{\Delta T}$.  The boundary condition at
Cartesian points near $\Gamma$ is enforced, with high-order accuracy,
as proposed in~\cite{appelo-peterson-2012}, on the basis of a certain
interpolation procedure which utilizes the given boundary values on
\(\Gamma\) as well as previously obtained solution values on interior
points of $H\cap \Omega$. The corresponding boundary condition at
Cartesian points near $\partial\Omega_{\Delta T}$, on the other hand,
are produced, for efficiency, by means of a special
procedure~\cite{fc_based_wave_solver_ff} which avoids a full
evaluation of the FC expansion for \(U_\mathrm{FC}\) at each
small-$\Delta t$ time interval and each boundary point, and which
constructs and use, instead, solutions similar to~\eqref{eq:aux_sol}
(without imposition of boundary conditions) but over some small square
regions contained in \(\Omega\) and centered at points on the interior
boundary of $\Omega\setminus\Omega_{\Delta T}$.

Combining the 2D-FC forwarded solution \(U_\mathrm{FC}\) and the near
boundary solution \(U_\mathrm{B}\) according to~\eqref{eq:ff_soln} the
desired FF numerical approximation of the solution \(u\) throughout
the Cartesian set \(H\cap \Omega\) is thus obtained up to time
$t=\Delta T$. Repeating this procedure as many times as necessary, the
solution can be advanced up to $t = n\Delta T$ for arbitrarily large
values of $n$, and, thus, up to an arbitrary final time $T$.

In view of the fact that auxiliary solutions
\(U_\mathrm{FC}(x, y, t)\) need to be computed only once per large
time step $\Delta T$, a significant improvement in the asymptotic
global computational cost per small time step $\Delta t$ results over
the cost required by classical finite-difference and finite-element
spatial discretizations.  Indeed, calling \(\delta\) the thickness of
the boundary region \(\Omega\setminus\Omega_{\Delta T}\), letting
\(n_c =[\delta/h]>0\) and assuming the Cartesian mesh $H\cap \Omega$
contains a total of \(\mathcal{O}(N)\) discretization points, it
follows that \(\Omega\setminus\Omega_{\Delta T}\) contains a total of
\(\mathcal{O}(\sqrt{N} n_c)\) grid points. As shown
in~\cite{fc_based_wave_solver_ff}, the optimum value of \(n_c\) is
\(\mathcal O(N^{1/4})\), so that the computational cost per time step
of the overall FF algorithm is \(\mathcal O(N^{3/4}\log N)\)
operations.  As shown in that contribution, further, owing to a
certain large multiplicative constant in front of the asymptotic cost
estimate for the time-stepping portion in the boundary region
\(\Omega\setminus\Omega_{\Delta T}\), large increases in $N$ are
necessary for the optimal $n_c$ value to increase by one or a few
units. Thus, in the numerical examples considered in the present
paper, for all of which we have $N\leq 4\cdot 10^6$, the value of
$n_c$ is set to a constant. This selection leads to an overall cost
estimate of approximately \(\mathcal O(\sqrt{N})\) operations for the
cases considered in this paper as the asymptotically large
\(\mathcal O(N\log N)\) FFT cost incurred by the algorithm has in fact
a limited impact in such cases.  A detailed discussion in this regard
is presented in~\cite{fc_based_wave_solver_ff}. The performance of the
resulting Fourier Forwarding method for a number of test cases is
demonstrated in~\Cref{sec:numer-results-ff} below.
\begin{figure}[h]
  \centering
  \captionsetup[subfigure]{labelformat=empty}
  \begin{subfigure}[t]{0.49\columnwidth}
    \centering
    \includegraphics[scale=0.52]{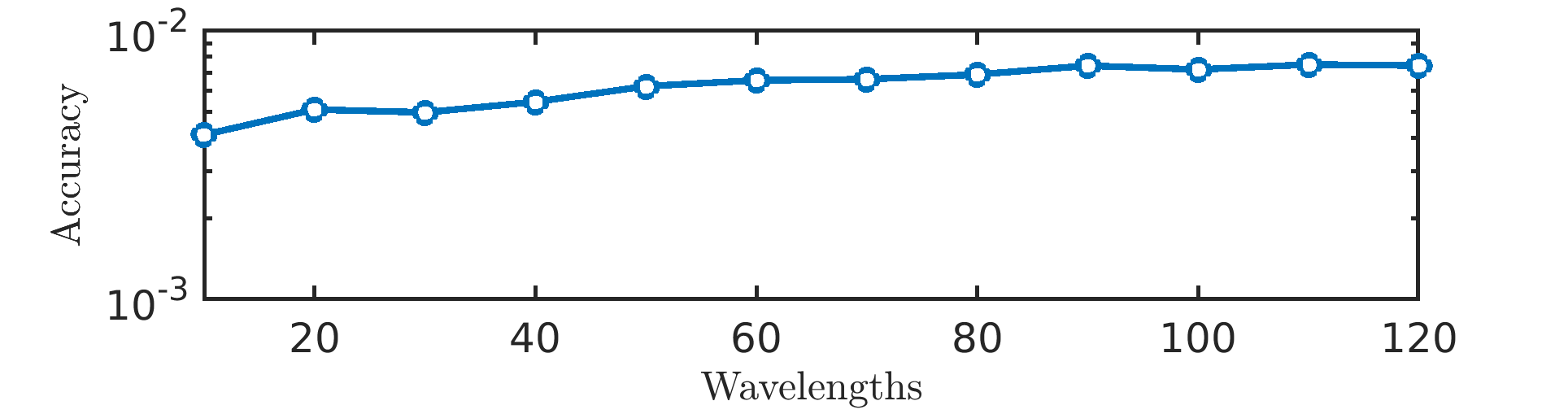}
  \end{subfigure}
  \begin{subfigure}[t]{0.49\columnwidth}
    \centering
    \includegraphics[scale=0.52]{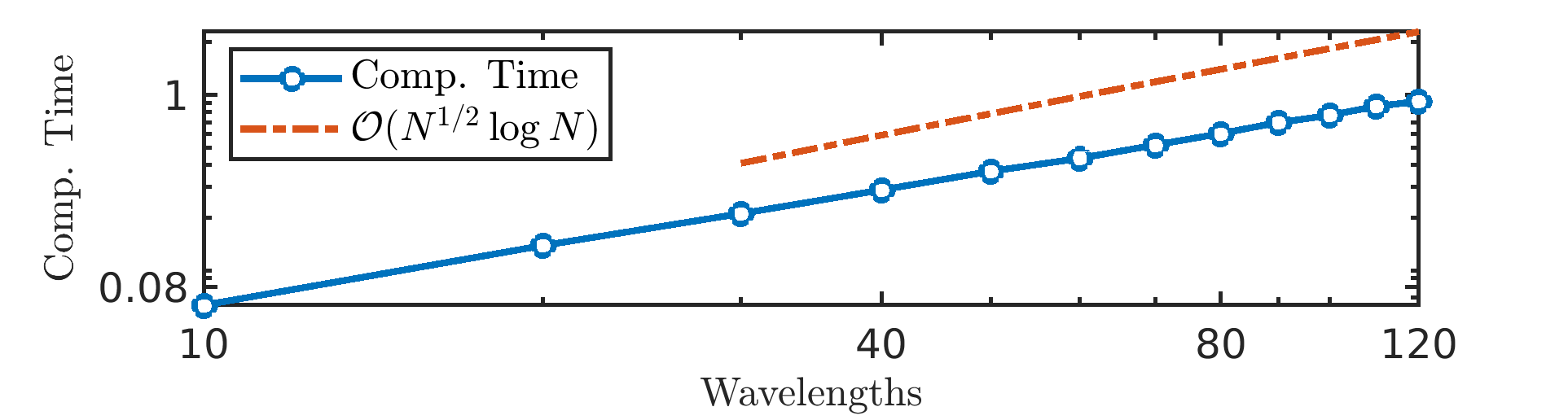}
  \end{subfigure}
  \caption{Accuracies and computing times (in seconds) for the problem
    considered in~\Cref{exm:ff_wave_eqn_2d_1}. Left graph: maximum
    absolute errors at time \(t=2\) for various wavelengths. Right
    graph: computing times required for each small-\(\Delta t\) time
    step.  Fifteen points per wavelength and \(n_c=6\) points across
    the boundary region were used for this numerical experiment. }~\label{fig:ff_2d_tsm2_15_ppw_nfd_6_wl_120}
\end{figure}

\subsubsection{Numerical Examples for the Fourier Forwarding method\label{sec:numer-results-ff}}
Two examples presented in this section demonstrate the character of
the FF method. The numerical examples presented in this section were
produced using Fourier continuations of order \(d=4\) for both the
1D-FC (for time stepping in the boundary region) and the 2D-FC (for
interior FC-forwarding) algorithms. In both the cases, the 2D-FC
parameter selections were made in accordance with \Cref{rmk:param_sel}
and the simulations were run in the computer described in
\Cref{rmk:comp_sys}. For time-stepping in the boundary region we have
utilized second order Taylor series method with the CFL number
\(0.125\).
\begin{example}[Accuracy and efficiency of the FF method]\label{exm:ff_wave_eqn_2d_1}
  In this example, we demonstrate the accuracy in the wave equation
  solution obtained via the FF method as well as other properties,
  namely, dispersionlessness and sublinear computing costs, enjoyed by
  the FF algorithm (\Cref{sec:fc_ff}) via the method of manufactured
  solutions. Here we consider the wave~\cref{eq:wave_2d} with the
  initial and (Dirichlet) boundary conditions taken such a way that
  the exact solution is given by
  \begin{align}
    \label{eq:ff_wave_eqn_2d_1}
    u(x, y, t) &= \cos(2\kappa (x + t) / 3) + \cos(\kappa (y + t)), 
  \end{align}
  on the unit disc
  $\Omega = \{(x,y)\in\mathbb{R}^2 : x^2 + y^2 \leq 1\}$. A fixed
  fifteen spatial points per wavelength, and a fixed number \(n_c=6\)
  of points across the boundary region (\(\delta = n_c h\)) have been
  utilized for this numerical experiment. In left graph of
  \Cref{fig:ff_2d_tsm2_15_ppw_nfd_6_wl_120}, we report the maximum
  absolute error (computed over all the Cartesian grid points within
  \(\Omega\)) of the solution, produced by the FF method, at time
  \(t=2\). The required computational times for each
  small-\(\Delta t\) time step for various spatial frequencies
  $\kappa$ (cf.~\eqref{eq:ff_wave_eqn_2d_1}) are reported in right
  graph of \Cref{fig:ff_2d_tsm2_15_ppw_nfd_6_wl_120}. As discussed at
  the end of \Cref{sec:fc_ff}, we have observed an
  \(\mathcal{O}(N^{1/2}\log{N})\) growth in the computational time as
  the size $N$ of the spatial discretization grows.
\end{example}
\begin{figure}[h]
  \captionsetup[subfigure]{labelformat=empty}
  \centering
    \begin{subfigure}[t]{0.25\columnwidth}
      \centering \includegraphics[width=0.8\textwidth]{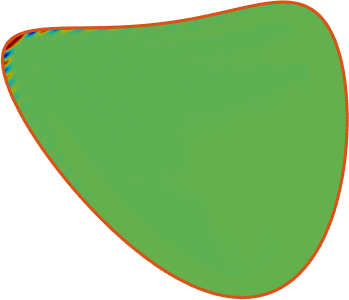} 
    \end{subfigure}
    \begin{subfigure}[t]{0.25\columnwidth}
      \centering \includegraphics[width=0.8\textwidth]{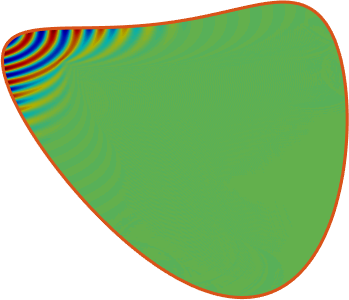}
    \end{subfigure}
    \begin{subfigure}[t]{0.25\columnwidth}
      \centering \includegraphics[width=0.8\textwidth]{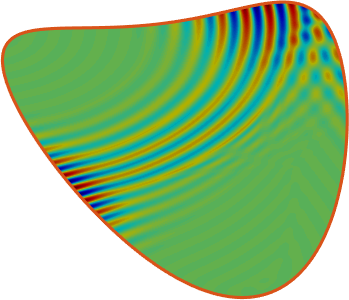}   
    \end{subfigure}
    \begin{subfigure}[t]{0.25\columnwidth}
      \centering \includegraphics[width=0.8\textwidth]{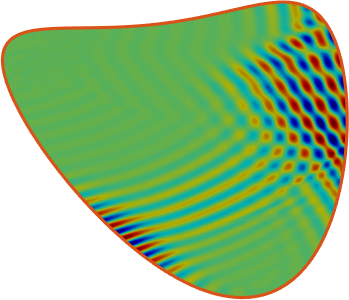}
    \end{subfigure}
    \begin{subfigure}[t]{0.25\columnwidth}
      \centering \includegraphics[width=0.8\textwidth]{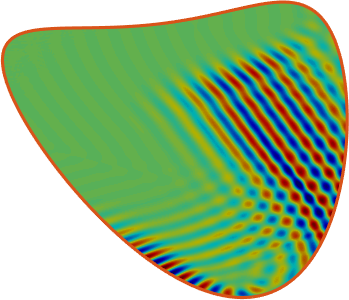}
    \end{subfigure}
    \begin{subfigure}[t]{0.25\columnwidth}
      \centering \includegraphics[width=0.8\textwidth]{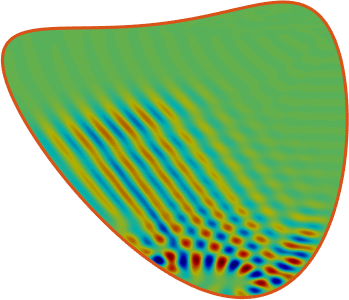}
    \end{subfigure}
    \caption{Fourier Forwarding method applied to an interior-domain
      wave propagation problem described
      in \Cref{exmp:ff_wave_eqn_2d_2}. The solution is shown, from top
      left to bottom right, after 100, 500, 2050, 2750, 3200 and 4300
      timesteps.}~\label{fig:kite_max_soln_value}
\end{figure}
\begin{example}[Interior-domain wave
  scattering]\label{exmp:ff_wave_eqn_2d_2}
  For an interior-domain graphical demonstration of the FF method, we
  consider the wave equation problem with boundary condition
  \begin{equation}\label{eq:bdry_cond_2}
    u = \cos(\kappa s) \exp( -s ^ 2 / \sigma ^ 2) g(t), 
  \end{equation}
  where
  \[s = (t -|\boldsymbol{r} - \boldsymbol{r_0}|) /|\boldsymbol{r} -
    \boldsymbol{r_0}|, \] and
  \begin{equation}\label{eq:windowing_function}
    g(t) =
    \begin{cases}
      0,& \textrm{for } t = 0,\\
      1 - \exp\left(2.0 \exp\left(\frac{-t_0 / t} {t / t_0 - 1.0}\right)\right),
      & \textrm{for } t < t_0,\\
      1,& \textrm{for } t \geqslant t_0,
    \end{cases}
  \end{equation}
  with \(\boldsymbol{r_0} = ( -1.1, -0.72)\), \(t_0 = 0. 05\) and
  \(\kappa = 20\pi\). The PDE domain is the interior of the
  kite-shaped curve considered in \Cref{exm:fc_2d_conv_kite}. A fixed
  number \(n_c=8\) of points across the boundary region and slightly
  over \(20\) points per wavelength have been utilized in this
  numerical experiment. Vanishing initial values of $u$ and $u_t$ at
  $t=0$ were used; note that both the imposed boundary
  values~\eqref{eq:bdry_cond_2} vanish for \(t \leq t_0 > 0\). For
  clear visibility, a version of the computed solution values at
  certain selected times scaled by the maximum value at that specific
  point in time are displayed in the left portion
  of~\Cref{fig:kite_max_soln_value}. (The scaling values range
  approximately between 1 and 3, and they are nearly equal to 1 for
  the first four images, all three in the upper row, and the leftmost
  image in the lower row.)
\end{example}
\section{Conclusions}\label{sec:conclusions}
This paper introduced a novel two-dimensional Fourier continuation
(2D-FC) method, for bi-periodic extension of functions defined on
arbitrary smooth two-dimensional domains. Applications to the Poisson
and wave-equation problem, including the development of the Fourier
Forwarding method, have resulted in numerical PDE solvers of high
orders of accuracy, and, most importantly, of extremely low numerical
dispersion.  Extensions of these methodologies to problems in higher
dimensions, and to problems on non-smooth domains, are left for future
work.

\section*{Acknowledgments}
This work was supported by NSF and DARPA through contracts DMS-1714169
and HR00111720035, and by the NSSEFF Vannevar Bush Fellowship under
contract number N00014-16-1-2808.


\end{document}